\documentclass[letterpaper, 10 pt, conference]{ieeeconf}  

\IEEEoverridecommandlockouts                              
\usepackage{graphicx}      
\usepackage{amsmath}
\usepackage{amssymb}
\usepackage{float}
\usepackage{graphicx}         
\usepackage{bm}
\usepackage{multirow}
\usepackage{subcaption}
\usepackage{algpseudocode}
\usepackage{algorithm}
\usepackage{makecell}
\usepackage{multirow}
\usepackage{comment}
\usepackage{breakcites}
\usepackage[table]{xcolor}
\let\labelindent\relax
\usepackage{enumitem}  
\usepackage{listings} 
\usepackage{xcolor}
\definecolor{mygreen}{RGB}{1, 70,20}

\newtheorem{thm}{Theorem}[section]
\newtheorem{lem}[thm]{Lemma}
\newtheorem{prop}[thm]{Proposition}

\newtheorem{rem}{Remark}

\title{\LARGE \bf
A Stabilizing NMPC Strategy \\ for a Class of Nonholonomic Systems with Drift
}

\author{Huu Thien Nguyen, Fernando A. C. C. Fontes, and Ionela Prodan
\thanks{This research is supported by  FCT/MCTES(PIDDAC), through projects 2020.07959.BD, 2022.02320.PTDC-KEFCODE, and 2022.02801.PTDC-UPWIND-ATOL (https://doi.org/10.54499/2022.02801.PTDC).}
\thanks{Huu Thien Nguyen and Fernando A. C. C. Fontes are with SYSTEC-ARISE, University of Porto, Porto, Portugal.
        {\tt\small huu-thien.nguyen@ieee.org, faf@fe.up.pt}}%
\thanks{Ionela Prodan is with Univ. Grenoble Alpes, Grenoble INP$^\dagger$, LCIS, F-26000, Valence, France.
        {\tt\small ionela.prodan@lcis.grenoble-inp.fr}. $^\dagger$ Institute of Engineering and Management Univ. Grenoble Alpes.}%
}

\begin{document}

\maketitle
\thispagestyle{empty}
\pagestyle{empty}

\begin{abstract}
In this paper, we present a \textcolor{black}{stabilizing} \textcolor{black}{Nonlinear Model Predictive Control (NMPC)} scheme tailored for a class of nonholonomic systems with drift, where the acceleration is laterally restrained. Examples include a mobile robot with drifting wheels on a planar surface or a spacecraft maneuvering in a vacuum.
\textcolor{black}{The novelty lies in the formulation of the terminal set, reachable from a significant distance from the equilibrium, and the terminal cost, represented as the integration of the stage cost.}
\textcolor{black}{
The proposed approach establishes essential steps for ensuring stability and feasibility guarantees. Simulation results substantiate the viability and effectiveness of the NMPC scheme.}
\end{abstract}
\section{INTRODUCTION}

Nonholonomic systems arise frequently in practice, everytime a mechanical device has a constraint on the motion in certain state-space directions \cite{kolmanovsky_DevelopmentsNonholonomic_1995}, examples of which are wheeled vehicles, fixed or rotary-wing aircraft, etc.
Prevailing research on nonholonomic systems predominantly revolves around the driftless unicycle model as evidenced by the survey \cite{nascimento_NonholonomicMobile_2018}.
Regarding a nonholonomic system with drift, \textcolor{black}{the most notable example is a} knife-edge moving on a planar plane \cite{bloch_ControllabilityStabilizability_1990}.
\textcolor{black}{A challenging property of nonholonomic systems is the uncontrollability of the linearized dynamics at the origin, which makes it impossible to stabilize them with a continuous feedback controller \cite{brockett_AsymptoticStability_1983}}.

\textcolor{black}{NMPC controllers, besides their established ability to manage constraints, are an adequate choice to deal with discontinuous state feedback \cite{a.c.c.fontes_DiscontinuousFeedbacks_2003}.}
Classical techniques to guarantee stability for an NMPC scheme are by using a terminal set constraint, defined as an ellipsoid \cite{chen_QuasiinfiniteHorizon_1998} or a polytope \cite{cannon_NonlinearModel_2003}, together with a terminal cost, defined as a Lyapunov function of a linearized system. However, these methods require linearization, which makes them an impossible choice to stabilize nonholonomic systems. 

Research on NMPC design for nonholonomic systems has delved into various applications, including 
i)  Differential-drive mobile robot and car-like: Worthmann et al. \cite{worthmann_ModelPredictive_2016}  formulate an NMPC for a unicycle mobile robot without terminal cost nor terminal constraint, where the stability is endorsed with a customized stage cost and a suitable prediction horizon. The differential-drive was also addressed with a stabilizing NMPC in \cite{a.c.c.fontes_DiscontinuousFeedbacks_2003,dongbinggu_StabilizingReceding_2005} and a robust NMPC in \cite{fontes_MinmaxModel_2003}. The car-like was addressed in \cite{fontes_GeneralFramework_2001}.
ii) Mobile robots with trailers: Rosenfelder et al. \cite{rosenfelder_ModelPredictive_2023} design NMPC schemes with kinematics-based non-quadratic costs for stabilizing driftless nonholonomic differential-drive vehicles.
iii)  Spacecraft attitude: Peterson et al. \cite{petersen_ModelPredictive_2017} construct an NMPC scheme for an underactuated spacecraft with two reaction wheels. Their NMPC scheme does not have a terminal cost, but they prove stability with the terminal state constraint equal to zero and the minimum prediction horizon of six steps. 
iii) Multicopters: 
\textcolor{black}{a multicopter is classified as a drifting nonholonomic system because, to transition between points in the planar plane or space, the multicopter needs to tilt to orient the thrust generated by its rotors. Subsequently, it tilts in the opposite direction to decelerate. In a specific example, as discussed in \cite{nguyen_NotesTerminal_2024}, an NMPC scheme was introduced with semi-globally asymptotic stability. However, it is important to highlight that this particular approach does not consider rotational dynamics (permitting instantaneous angle changes), and instead, a feedback linearization controller is employed to linearize the multicopter around the equilibrium.}

\textcolor{black}{To the best of the authors' knowledge, this is the first work to design a stabilizing NMPC for a class of drift nonholonomic systems with no lateral acceleration.}

\textbf{Notation}:
Let $\bm{0}$ be the zero matrix or vector of appropriate dimension depending on the context,  
$\mathbb{R}^{+}$ be the set of positive real numbers.
We denote $\mathbf{x}(t; t_0, \mathbf{x}_{t_0} ,\mathbf{u})$ the state $\mathbf{x}$ at time $t$, starting from $ \mathbf{x} = \mathbf{x}_{t_0}$ at $t=t_0$, with the control law $\mathbf{u}$ for $t \geq t_0$.
\textcolor{black}{We also use $x_{t_i} \triangleq x(t_i)$ to present the value of the function $x(t)$ at the time index $t_i$.}
 \textcolor{black}{The outline of the paper is as follows:
 we first present a spacecraft model, which is a drifting nonholonomic system in Section \ref{sec_model}.
 Next, we introduce a terminal set, together with its associated auxiliary controller in Section \ref{sec_auxiliary_controller}.
 Section \ref{sec_NMPC_framework} establishes the NMPC scheme and the main results of the paper, while we demonstrate them through simulations in Section \ref{sec_simulation}. 
 Finally, we state the conclusions and suggest some future work in Section \ref{sec_conclusion}.
}
\section{SPACECRAFT MODEL}
\label{sec_model}
\begin{figure}[ht!]
	\centering
	\includegraphics[width=0.8\linewidth]{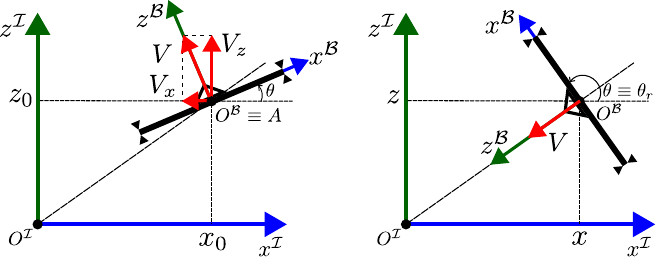}
	\caption{Spacecraft with an arbitrary pitch angle $\theta$ (left) and when its velocity vector $\bm{V}$ and its vertical axis $x^\mathcal{B}$ point to the origin (right).}
	\label{fig_Spacecraft_2D_Full_Model}
\end{figure}

Consider a spacecraft operating in vacuum (with no air resistance and no gravity) as in Fig. \ref{fig_Spacecraft_2D_Full_Model}.
The spacecraft is constrained to a plane where it can only move forward, move backward, and rotate around its center of mass by using two bi-directional thrusters. The same model can also represent a mobile robot with drifting wheels operating on a slippery plane. We will describe the dynamics of this system and show that it is a nonholonomic system with drift.  

We define two coordinate frames, the fixed inertial frame $\mathcal{I}: \{O^\mathcal{I},x^\mathcal{I}, z^\mathcal{I} \}$ with basis $(\bm{e}_x^\mathcal{I},\bm{e}_z^\mathcal{I})$, and the body frame $\mathcal{B}:\{O^\mathcal{B},x^\mathcal{B}, z^\mathcal{B} \}$ whose $O^\mathcal{B}$ coincides with the center of mass of the spacecraft, and with the basis $(\bm{e}_x^\mathcal{B},\bm{e}_z^\mathcal{B})$.

\vspace{-0.3cm}


The state-space model of the aircraft is:
\begin{align}
\label{eqn_SC_state_space_dynamics}
	\begin{cases}
		\dot{x}(t) &= V_x(t),  \\
		\dot{z}(t) &= V_z(t),  \\
		\dot{\theta}(t) &= \omega(t), \\
		\dot{V}_x(t) &= -a(t) \sin \theta(t),  \\
		\dot{V}_z(t) &= a(t) \cos \theta(t).  \\
	\end{cases}
\end{align}
Here, 
 the state and control vectors are
$\mathbf{x} \triangleq [x,z,\theta, V_x, V_z]^\top$, $\mathbf{u} \triangleq [a,\omega]^\top$, \textcolor{black}
{where $x,z$ are the position of $O^\mathcal{B}$ in $\mathcal{I}$, $V_x$ and $V_z$ are its velocity components in the same inertial frame.}
We also consider the pitch angle between the body frame ${x^\mathcal{B}}$ axis and the inertial frame ${x^\mathcal{I}}$ axis, i.e., $\theta(t) = \angle(\bm{e}_x^\mathcal{B}, \bm{e}_x^\mathcal{I})$, $\theta(t) \in [-\pi,\pi]$.
The control components,
 $\omega$ and  $a$, are the angular speed and the translational acceleration, respectively, with
$a(t) \in [-a_m,a_m]$, $\omega(t) \in [-\omega_m,\omega_m]$. 
%
%
The input constraint set is then: 
\begin{equation}
	\label{eqn_SC_dynamics_input_constaint_set}
\mathcal{U} \triangleq \{\mathbf{u}=\in \mathbb{R}^2 
\mid
-[a_m,\omega_m]^\top \leq \mathbf{u} \leq [a_m,\omega_m]^\top\}.
\end{equation}
We define the \emph{displacement vector} $\bm{r} \triangleq x(t) \bm{e}_x^\mathcal{I} + z(t) \bm{e}_z^\mathcal{I} $, with $r(t) \triangleq \sqrt{x(t)^2+z(t)^2}$ being the distance between the spacecraft and the origin. Similarly, we define the \emph{velocity vector} $\bm{V} \triangleq V_x(t) \bm{e}_x^\mathcal{I} + V_z(t) \bm{e}_z^\mathcal{I} $, with $V(t) \triangleq \sqrt{V_x(t)^2+V_z(t)^2}$.

The dynamics \eqref{eqn_SC_state_space_dynamics} are an example of a nonholonomic system with drift \cite{godhavn_SteeringClass_1999}, which could be written in a compact form as:
\begin{equation}
	\label{eqn_dynamics_nonholonomics}
	{\mathbf{\dot{x}}}(t) 
	= f(\mathbf{x}(t),\mathbf{u}(t)) 
	= f_1(\mathbf{x}(t)) 
	+ \sum_{i=1}^{2}
	g_i(\mathbf{x}(t))u_i(t),
\end{equation}
where the drift term is $f_1(\mathbf{x}(t)) = [V_x(t), V_z(t),0, 0, 0]^\top$, while $g_1(\mathbf{x}(t)) = [0, 0, 0, \sin\theta(t), \cos\theta(t)]^\top$, $g_2(\mathbf{x}(t)) = [0, 0, 1, 0, 0]^\top$, $u_1(t) = a(t)$, and $u_2(t) = \omega(t)$. 
\textcolor{black}{Moreover, the} dynamics \eqref{eqn_SC_state_space_dynamics} admit a nonholonomic constraint:
$
	\label{eqn_SC_state_space_dynamics_nonholonomic_constraint}
	\dot{V}_x(t) \cos\theta(t) + \dot{V}_z(t) \sin \theta(t) = 0,
$
which prevents the spacecraft from accelerating in the $x^\mathcal{B}$ direction. 
As is the case with any nonholonomic system, the linearization of the dynamics \eqref{eqn_SC_state_space_dynamics} at $(\mathbf{x},\mathbf{u}) = (\bm{0},\bm{0})$ is neither controllable nor stabilizable, which forbids the use of any methods based on linearization, including dual-mode \cite{chen_QuasiinfiniteHorizon_1998} and a vast number of MPC works based on it. 
\section{AUXILIARY CONTROL STRATEGIES}
\label{sec_auxiliary_controller}
In this section, we propose a stabilizing strategy for the nominal model. The strategy consists of first driving the state to a set containing the space directions where the nonholonomic systems are easier to stabilize, and then moving to the origin by sliding along this set.
Such a process has similarities with sliding mode control, where we first drive the system to a sliding set and then drive the system to the origin along that set. The auxiliary strategy will not necessarily be applied to the system, but will only be used to show that an NMPC scheme is stabilizing. The set mentioned will be the terminal set in the NMPC scheme. This set can be defined to be the set of states with a heading angle pointing towards the origin (see the right side of Fig. \ref{fig_Spacecraft_2D_Full_Model}), from where we can easily reach the origin simply by moving straight on.

Starting from a point A $(x_0,z_0)$ stopped, with $\mathbf{x}_A = [x_0,z_0,0,0,0]^\top$, the goal is to move the spacecraft to the equilibrium point $\mathbf{x} = \bm{0}$.
A possible auxiliary stabilizing strategy (SS) is: 
\begin{enumerate}[label={SS\arabic*},leftmargin=21pt]
	\item \label{eqn_SS_1} Rotate the spacecraft until its vertical axis \textcolor{black}{ $\bm{z^\mathcal{B}}$ points to the origin (see Fig. \ref{fig_Spacecraft_2D_Full_Model})}.	
	\item \label{eqn_SS_2} Move forward the target with an acceleration profile $a_\text{ref}(t)$ (and consequently with a velocity profile $v_\text{ref}(t)$) until reaching the origin of the plane $(x,z) = (0,0)$.
	\item \label{eqn_SS_3} Rotate back to the equilibrium pitch angle $\theta = 0$ and then stop.
\end{enumerate}
To implement this strategy it is convenient to define the \emph{reference pitch angle}, which for any position  $(x,z)$ gives the pitch angle that points the thrust to the origin:
\begin{equation}
	\label{eqn_theta_ref}
	\theta_r(x,z) =
	\begin{cases}
		0 & \mbox{if\ } x=0,z\leq 0;\\
		\pi  &\mbox{if\ } x=0,z > 0;\\
		\tan^{-1}(-x/z)  &\mbox{if\ } x \neq 0, z <0 ;\\
		\pi - \tan^{-1}(x/z)   &\mbox{if\ } x\neq0, z >0;\\
	\end{cases}
\end{equation}
\textcolor{black}{Let us define the} \emph{auxiliary controller} that implements the \textcolor{black}{above} strategy:
\begin{align}
		\label{eqn_SC_auxiliary_controller}
&	\mathbf{u}_\text{aux}    = [a, \omega]^\top = \nonumber  \\
&	\left\{
	\begin{array}{lll}
		a = 0,                 &\omega = -\omega_m &\text{ if } \theta  > \theta_r;    \\
		a = 0,                 &\omega = \omega_m  &\text{ if } \theta  < \theta_r;  \\
		a = a_{r}(t), &\omega =0          &\text{ if }  \theta = \theta_r, (r,V)\neq \bm{0};    \\
		a=0,                   &\omega =0          &\text{ if } \mathbf{x}=\bm{0}; \\
	\end{array}\right.
\end{align}
which is equivalent to \emph{Rotate clockwise}, \emph{Rotate anticlockwise}, \emph{Move forward}, and \emph{Stop}. The rotation when $\theta \neq \theta_r$  could be written in the shortened form with the help of the sign function: 
$
\mathbf{u}_\text{aux}
=[a, \omega]^\top 
= [0, \text{sign}(\theta_r-\theta) w_m]^\top   \text{ if }   \theta \neq \theta_r.  
$
 It is clear that the aforementioned control strategy satisfies the input constraint \eqref{eqn_SC_dynamics_input_constaint_set}. Moreover, it is a discontinuous control law, as expected by Brockett's theorem \cite{brockett_AsymptoticStability_1983}.

\begin{figure}[ht!]	
	\centering
	\begin{subfigure}[b]{0.49\linewidth}
		\centering
		\includegraphics[width=\textwidth]{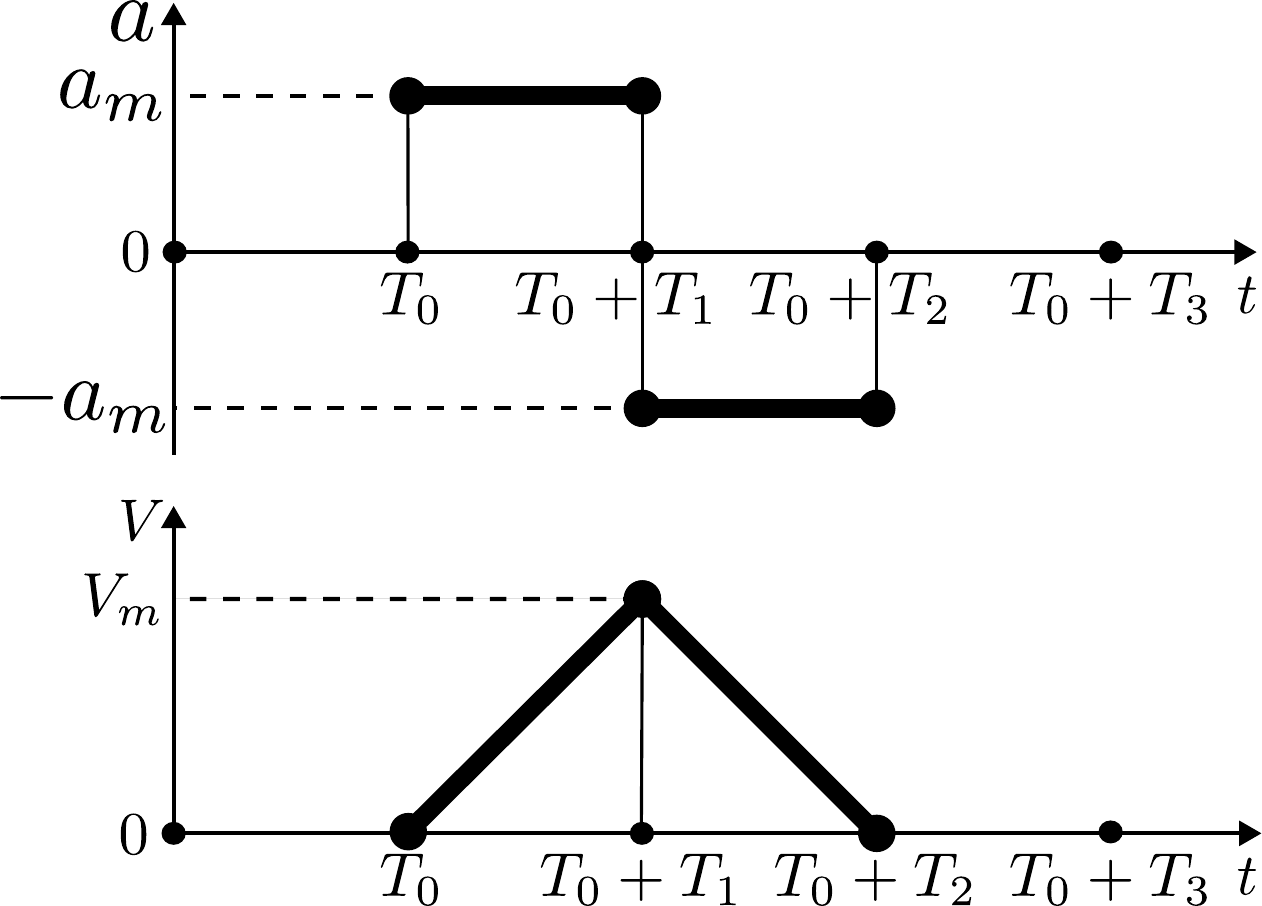}
		\caption{$a$ and $V$}
		\label{fig_SC_fwd_acc_vel}
	\end{subfigure}
	\hfill
	\begin{subfigure}[b]{0.49\linewidth}
		\centering
		\includegraphics[width=\textwidth]{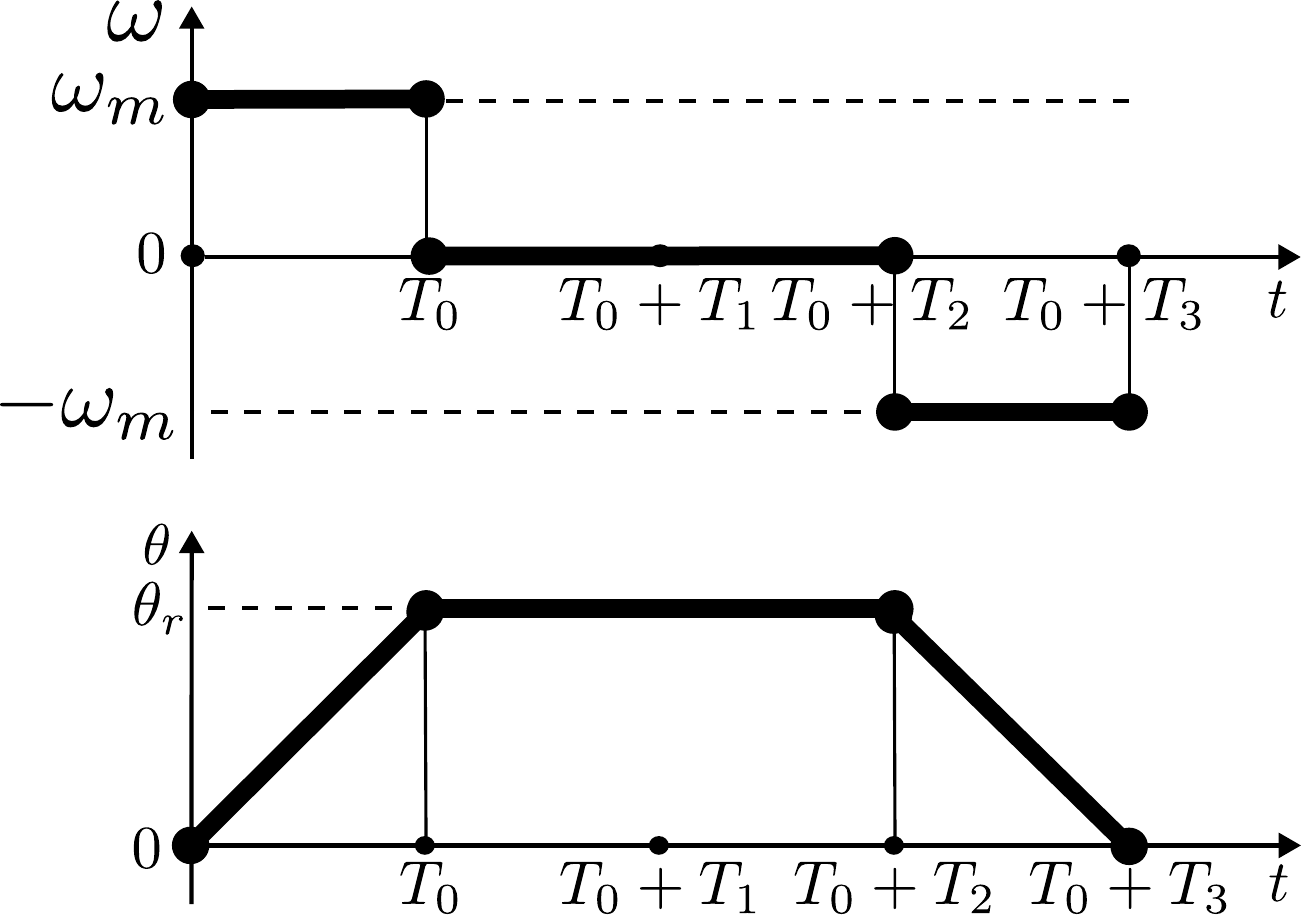}
		\caption{$\omega$ and $\theta$}
		\label{fig_SC_fwd_omega_theta}
	\end{subfigure}
	\caption{\textcolor{black}{Proposed} auxiliary stabilizing controller and states.} 
	\label{fig_SC_fwd_acc_vel_omega_theta}
\end{figure}

The acceleration profile $a_r(t)$ when starting from a stopped position $(V_x,V_z)=(0,0)$ is (see Fig. \ref{fig_SC_fwd_acc_vel}):
\begin{align*}
	a_r(t) =  
	\left\{
	\begin{array}{lll}
		a_m  & \hspace{-5pt} \text{if } t \in [T_0, T_0 + T_1];  &\text{(Accelerate)} \\
		-a_m & \hspace{-5pt}  \text{if } t \in [T_0 + T_1, T_0 + T_2]   &\text{(Decelerate)} 
	\end{array}\right.
\end{align*}

One of the three cases for the auxiliary stabilizing controller \eqref{eqn_SC_auxiliary_controller} is presented in Fig. \ref{fig_SC_fwd_acc_vel_omega_theta}: 
i) for $t \in [0,T_0+T_2]$, the reference pitch angle $\theta_r$ is positive.
ii) for $t \in [T_0+T_2,T_0+T_3]$, when the spacecraft is at the origin, i.e. $(x,z)=(0,0)$, the reference pitch angle $\theta_r$ is 0.

We want the spacecraft to accelerate then decelerate with equal intervals (i.e. $T_2=2T_1$), and it travels the distance $r_0$ from $T_0$ to $T_0 + T_2$, \textcolor{black}{following the relation} $r_0/2 = a_m T_1^2/2$ (see Fig. \ref{fig_SC_fwd_acc_vel}), \textcolor{black}{which leads to:}
\begin{equation}
	\label{eqn_SC_T1_T2_explicit_main}
	T_1 =  \sqrt{\dfrac{r_0}{a_m}},  \quad
	T_2 =  2\sqrt{\dfrac{r_0}{a_m}}.
\end{equation}
Regarding the rotations, the rotation from $\theta=0$ to $\theta=\theta_r$ for $t \in [0,T_0]$ and the rotation from $\theta=\theta_r$ back to $\theta=0$ for $t \in [T_0+T_2,T_0+T_3]$ give:
\begin{equation}
	T_0 = \dfrac{|\theta_r|}{ w_m}, \quad 
	T_3 = 2\sqrt{\dfrac{r_0}{a_m}} + \dfrac{|\theta_r|}{ w_m}.
	\end{equation}
Noting that the spacecraft is symmetrical in the longitudinal axis $\bm{x^\mathcal{B}}$, we can similarly define a second stabilizing strategy (SS') where the vertical axis $\bm{z^\mathcal{B}}$ points away from (rather than pointing to) the origin and then the spacecraft moves backward to the origin. The details are in Appendix \ref{appendix_second_auxiliary_controller}.

\textcolor{black}{Next, let us} define the \emph{terminal set} to be the set of states for which the spacecraft's vertical axis 
\textcolor{black}{ $\bm{z^\mathcal{B}}$ points toward (or points away from) the origin}
and the spacecraft is moving toward (or reversing backward) the origin:
	\begin{align}
			\label{eqn_SC_NMPC_terminal_set_def}
	\mathcal{X}_{f} \triangleq   \{(x,z,\theta,V_x,V_z) \in \mathbb{R}^2 \times [-\pi,\pi]  \times  \mathbb{R}^2 \mid &\nonumber \\
		\left.	\bm{r} = -\kappa_1 \bm{V} = \kappa_2  \bm{e}_z^\mathcal{B}, 
		\kappa_1 \in R^{+}, \kappa_2 \neq 0, V^2 \leq 2a_m r
		\right\}.  &
	\end{align}

\textcolor{black}{This set $\mathcal{X}_f$ will serve as terminal region in the NMPC framework. 
The last condition in the definition of the terminal set imposes that the velocity is sufficiently small so that the spacecraft has room to decelerate and stop at the origin of the plane $(x,z) = (0,0)$. This is further discussed in Section \ref{subsec_Stability_Analysis}.}

 We can see that with this set definition, the auxiliary strategy \eqref{eqn_SC_auxiliary_controller} makes $\mathcal{X}_{f}$ invariant 
\textcolor{black}{with $\kappa_2 <0$}, while the auxiliary strategy SS' makes $\mathcal{X}_{f}$ invariant \textcolor{black}{with $\kappa_2 > 0$}. 
\textcolor{black}{In this paper, we focus on the analysis of the auxiliary strategy SS, we also obtain similar results with the auxiliary strategy SS'}.
The auxiliary strategy SS is equivalently defined as:
\begin{enumerate}[label={AC\arabic*},leftmargin=24pt]
	\item Rotate to $\mathcal{X}_{f}$, i.e. rotate to align the pitch angle $\theta = \theta_r(x,z)$, with the angular rate $\omega = \text{sign}(\theta_r-\theta) w_m$.
	\item \label{AC_2} Move straight to $(x,z)=(0,0)$ within $\mathcal{X}_f $ following the acceleration profile: $a=a_m \text{ for } t \in [t_f,t_f+t_1]$ and $a=-a_m \text{ for } t \in [t_f+t_1,t_f+t_2]$.
	\item \label{AC_3} Rotate to $(x,z,\theta)=(0,0,0)$ with the angular rate  $\omega = \text{sign}(\theta_r-\theta) w_m = \text{sign}(-\theta) w_m = -\text{sign}(\theta) w_m$ $\text{ for } t \in [t_f+t_2,t_f+t_3]$ and stop.
\end{enumerate}
The values of $t_1$, $t_2$, and $t_3$ are computed in Section \ref{subsec_explicit_terminal_cost}. 

\begin{figure}[ht!]
	\centering
	\includegraphics[width=0.4\linewidth]{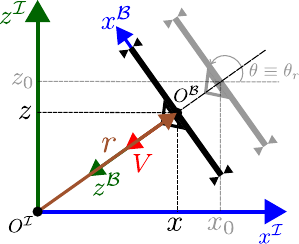}
	\caption{The spacecraft when it is inside the terminal set $\mathcal{X}_f$ and points to the origin.}
	\label{fig_Spacecraft_2D_Full_Model_r_V}
\end{figure}

Figure \ref{fig_Spacecraft_2D_Full_Model_r_V} shows a snapshot of the spacecraft when it ``slides forward'' inside the terminal set \textcolor{black}{toward the origin}, with the velocity vector $\bm{V}$ in red, and the displacement vector $\bm{r}$ in brown. Within the terminal set $\mathcal{X}_f$, we obtain simpler dynamics for the spacecraft \textcolor{black}{defined in \eqref{eqn_SC_state_space_dynamics}}:
\begin{subequations}
	\label{eqn_SC_terminal_set_dynamics}
	\begin{align}
		\dot{r}(t) &= -V(t),  \\
		\dot{\theta}(t) &= \omega(t), \\
		\dot{V}(t) &=  a(t), 
	\end{align}
\end{subequations}
or, in the \textcolor{black}{compact} form 
with $\mathbf{x}_{ts}(t) = [r(t) , \theta(t), V(t)]^\top$ and $\mathbf{u}_\text{aux}(t)$ as in \eqref{eqn_SC_auxiliary_controller}:
\begin{equation}	
  \label{eqn_SC_terminal_set_dynamics_shortened}
	\dot{\mathbf{x}}_{ts}(t) = f_{ts}(\mathbf{x}_{ts}(t),\mathbf{u}_\text{aux}(t)).
\end{equation}
\section{NMPC FRAMEWORK}
\label{sec_NMPC_framework}
At the time instant $t$, with the state $\mathbf{x}_{t}$, we solve an open-loop optimal control problem (OCP) :
\begin{equation}
	\label{eqn_SC_NMPC_problem_formulation}
	\underset{{\bar{\mathbf{u}}(\cdot)}}{\text{min}}   \ J 
	\left(\bar{\mathbf{x}} ,\bar{\mathbf{u}} \right)
	= \hspace{-2pt}  \int_{t}^{t+T_p}  \hspace{-5pt}
	L (\bar{\mathbf{x}}(\tau),\bar{\mathbf{u}}(\tau) )  d \tau
	+ F(\bar{\mathbf{x}}(t+T_p))   ,
\end{equation}
\vspace{-14pt}
\begin{subequations}
	\label{eqn_SC_NMPC_problem_formulation_constraints}
	\begin{align}       
		\text{subject to: }		 
		& \dot{\bar{\mathbf{x}}}(\tau) 
		=
		f(\bar{\mathbf{x}}(\tau),\bar{\mathbf{u}}(\tau)), \ \tau \in [t,t+T_p]  \\
		\label{eqn_SC_NMPC_problem_formulation_state_input_constraints}
		&\bar{\mathbf{u}}(\tau)  \in \mathcal{U}, \ \tau \in [t,t+T_p] \\
		\label{eqn_SC_NMPC_problem_formulation_initial_states}
		& \bar{\mathbf{x}}(t) = \mathbf{x}_{t}, \\
		& \bar{\mathbf{x}}(t+T_p) \in \mathcal{X}_f.
	\end{align}
\end{subequations}
Here $\bar{\mathbf{x}}(\tau) = \mathbf{x}(\tau; t, \mathbf{x}_{t} ,\bar{\mathbf{u}})$ is the predicted state at the time $\tau$, starting from $\mathbf{x}_{t}$ at $t$, under the predicted control input $\bar{\mathbf{u}}$,
which is the optimal solution to the OCP.
The set $\mathcal{U}$ is input constraint set in  \eqref{eqn_SC_dynamics_input_constaint_set}, 
 $\mathcal{X}_f$ is terminal set defined in \eqref{eqn_SC_NMPC_terminal_set_def}, and $T_p$ is the prediction horizon.
The first part of the solution, $\bar{\mathbf{u}}(\tau), \tau \in [t,t+\delta ]$, is applied to the system, then the time shifts $t \leftarrow t+\delta$, and the problem \eqref{eqn_SC_NMPC_problem_formulation}--\eqref{eqn_SC_NMPC_problem_formulation_constraints} is solved repeatedly.
The MPC law $\mathbf{u}^\star(\cdot)$ \textcolor{black}{is constructed }by the concatenation of the solutions of \eqref{eqn_SC_NMPC_problem_formulation}--\eqref{eqn_SC_NMPC_problem_formulation_constraints} that are applied to the system \cite{fontes_GeneralFramework_2001,grune_NonlinearModel_2017}.

We define respectively the stage cost $L$ and the terminal cost $F$ in \eqref{eqn_SC_NMPC_problem_formulation} as follows:
\begin{subequations}
\begin{align}
    	\label{eqn_SC_NMPC_problem_formulation_stage_cost}
	L({\mathbf{x}}_\tau,{\mathbf{u}}_\tau )
	&=   x_\tau^2 + z_\tau^2 + {V_x}_\tau^2 + {V_z}_\tau^2 + \theta_\tau^2, \\
 	\label{eqn_SC_NMPC_problem_formulation_terminal_cost}
	F({\mathbf{x}}_{t_f} ))  
	&=  \int_{t_f}^{\infty} 
	L (\bar{\mathbf{x}}(\tau;t_f,{\mathbf{x}}_{t_f},\mathbf{u}_\text{aux} ),
	 \mathbf{u}_\text{aux}(\tau)) d\tau.
\end{align}
\end{subequations}
The terminal cost defined this way is guaranteed to be finite and can be explicitly computed since the auxiliary controller drives the state to the origin in finite time. 

We provide in Algorithm \ref{algorithm_NMPC_SC_design} the steps to design an NMPC controller with guaranteed stability.
\begin{algorithm}
	\caption{NMPC for a drifting nonholonomic vehicle in the planar plane with laterally restricted acceleration}
	\label{algorithm_NMPC_SC_design}
	\begin{algorithmic}[1]
		\Statex \textbf{data} The dynamics $f$ in \eqref{eqn_dynamics_nonholonomics}, the initial state $\mathbf{x}_0$, the acceleration and angular velocity limits $a_m$, $\omega_m$, the error tolerance $\epsilon_r$.
  \State Define the MPC design parameters ($L,F,T_p,\mathcal{X}_f$) satisfying the stability conditions \ref{SC_1}--\ref{SC_Lyapunov} defined below.
		\State Set the time $t=0$.
		\While {$r^2 \geq \epsilon_r$} \label{eqn_epsilon_r} 
		\State Obtain the current state $\mathbf{x}_t$.
		\State Solve the OCP \eqref{eqn_SC_NMPC_problem_formulation}--\eqref{eqn_SC_NMPC_problem_formulation_constraints}. 
		\State  Assign the first segment of the optimal input, \textcolor{black}{$\bar{\mathbf{u}}(\tau), \tau \in [t,t+\delta]$ } to  $\mathbf{u}^\star(\cdot)$ and apply it to the plant.
		\State Shift the time $t \leftarrow  t+\delta$.
		\EndWhile
		\Statex \textbf{result} NMPC law $\mathbf{u}^\star(\cdot)$
	\end{algorithmic}
\end{algorithm}
\subsection{Explicit expression for the terminal cost}
\label{subsec_explicit_terminal_cost}
\begin{figure}[ht!]
	\centering
	\begin{subfigure}[b]{0.49\linewidth}
		\centering
		\includegraphics[width=\textwidth]{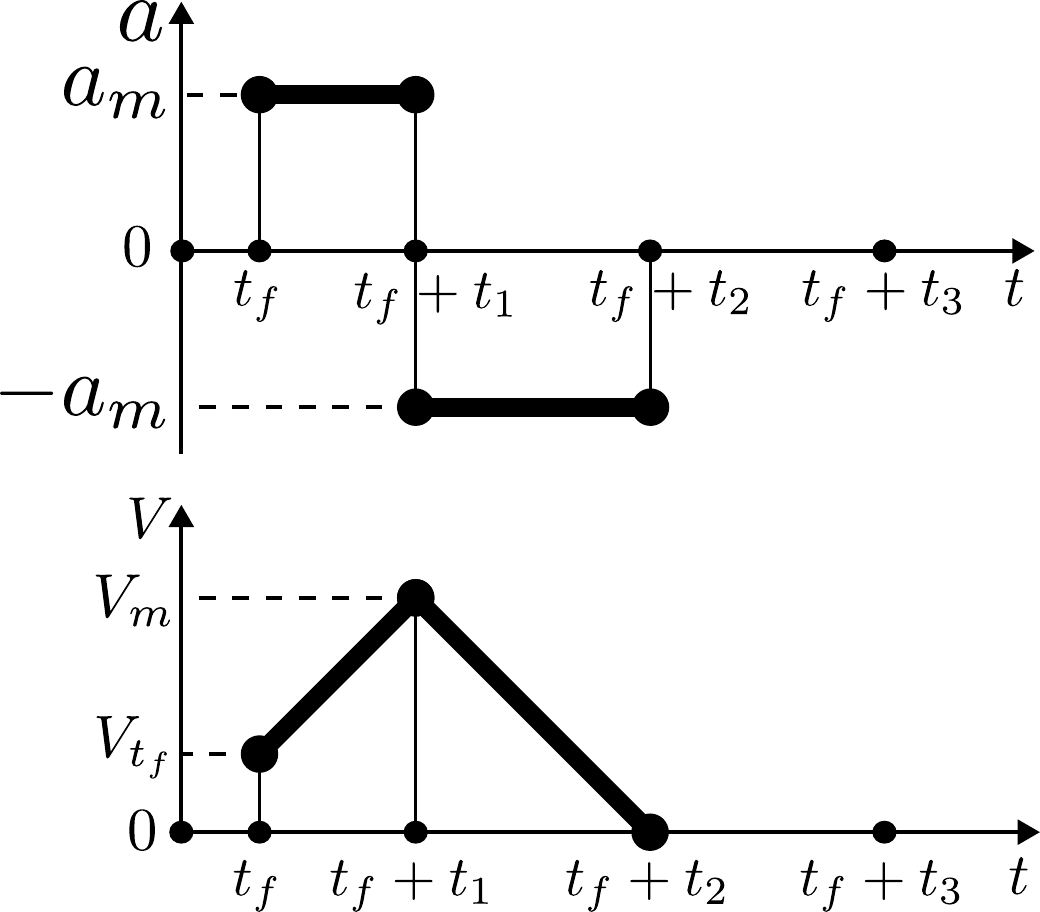}
		\caption{$a$ and $V$}
		\label{fig_SC_terminal_set_acc_vel}
	\end{subfigure}
	\hfill
	\begin{subfigure}[b]{0.49\linewidth}
		\centering
		\includegraphics[width=\textwidth]{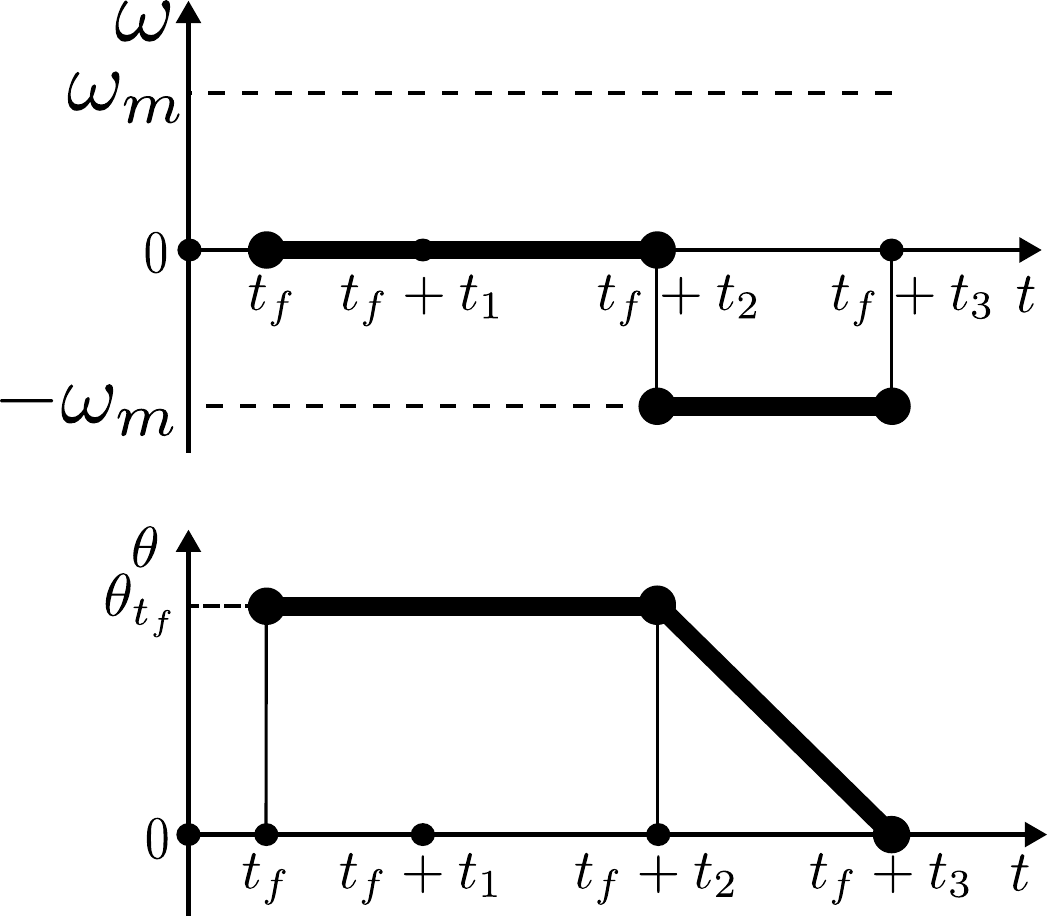}
		\caption{$\omega$ and $\theta$}
		\label{fig_SC_terminal_set_omega_theta}
	\end{subfigure}
	\caption{The auxiliary stabilizing controller AC and states inside the terminal set.}
	\label{fig_SC_terminal_set_acc_vel_omega_theta}
\end{figure}

Figure \ref{fig_SC_terminal_set_acc_vel_omega_theta} presents the concept of the auxiliary controller inside an NMPC iteration when $\theta_{t_f} > 0 $. After the time stamp $t_f$, the auxiliary controller $\mathbf{u}_\text{aux}$, which is never applied to the system, is used to compute the terminal cost $F$ and verify the stability conditions. 
For $t \in [t_f,t_f+t_1]$, $\mathbf{u}_\text{aux} =[a_m,0]^\top$ to accelerate the spacecraft to its maximum velocity $V_m$, 
then for  $t \in [t_f+t_1,t_f+t_2]$, $\mathbf{u}_\text{aux} = [-a_m,0]\top$ to decelerate the spacecraft to zero velocity. 
Finally, for $t \in [t_f+t_2,t_f+t_3]$,  $\mathbf{u}_\text{aux} 
= [0,\text{sign}(\theta_r-\theta_{t_f})\omega_m]\top 
=[0, \text{sign}(-\theta_{t_f})\omega_m]\top
=[0, -\text{sign}(\theta_{t_f})\omega_m]\top$ to rotate the spacecraft's pitch angle from $\theta_{t_f}$ to $0$. The time stamps $t_1$, $t_2$, and $t_3$ are chosen such that, the velocity $V$ and the distance $r$ at $t_f+t_2$ are 0, while from $t_f+t_2$ to $t_f+t_3$, the pitch angle goes from $\theta_{t_f}$ to $0$. 
\begin{lem}
	\label{lemma_t1_t2_t3}
	The time stamps $t_1$, $t_2$, and $t_3$ in the auxiliary strategy \ref{AC_2}--\ref{AC_3} when starting from $(r,V)=(r_{t_f}, V_{t_f})$ are calculated as follows:
 \begingroup
\allowdisplaybreaks
\begin{subequations}
		\label{eqn_SC_t1_t2_explicit_main}
		\begin{align}
			\label{eqn_SC_t1_explicit_main}
			t_1 &= -\dfrac{V_{t_f}}{a_m} + \dfrac{\sqrt{V_{t_f}^2 + 2a_m r_{t_f}}}{a_m \sqrt{2}}, \\	
			\label{eqn_SC_t2_explicit_main}
			t_2 &= -\dfrac{V_{t_f}}{a_m} + \dfrac{\sqrt{2} \sqrt{V_{t_f}^2 + 2a_m r_{t_f}}}{a_m },\\
						\label{eqn_SC_t3_explicit_main}
			t_3 &=   \frac{ |\theta_{t_f}|}{ \omega_m} -\dfrac{V_{t_f}}{a_m} + \dfrac{\sqrt{2} \sqrt{V_{t_f}^2 + 2a_m r_{t_f}}}{a_m }.
		\end{align}
	\end{subequations}
 \endgroup
	\vspace{0.05pt}
\end{lem}
\proof The proof is in Appendix \ref{appendix_terminal_cost_steps}. \hfill\ensuremath{\square}
\textcolor{black}{\begin{rem}
	The time stamps $t_2$ and $t_3$ are always non-negative. Hence, the time stamp $t_1$ is non-negative when: 
	\begin{equation}
V_{x_{t_f}}^2 +V_{z_{t_f}}^2 \leq 2 a_m \sqrt{x_{t_f}^2+z_{t_f}^2},
	\end{equation}
	which explains the \textcolor{black}{last condition} in \eqref{eqn_SC_NMPC_terminal_set_def}.
\end{rem}}

Inside the terminal set, with the state variables $\mathbf{x}_{ts}(t) = [r(t) , \theta(t), V(t)]^\top$ in \eqref{eqn_SC_terminal_set_dynamics_shortened}, we write the stage cost \eqref{eqn_SC_NMPC_problem_formulation_stage_cost} as:
\begin{align}
  L(\bar{\mathbf{x}}(\tau),\bar{\mathbf{u}}(\tau) )
&=   x(\tau)^2 + z(\tau)^2 + V_x(\tau)^2 + V_z(\tau)^2 + \theta(\tau)^2  \nonumber \\ 
& = r(\tau)^2 + V(\tau)^2 + \theta(\tau)^2.
\end{align}
Using this new stage cost, with the auxiliary controller \ref{AC_2}--\ref{AC_3}, we derive the exact formula for the terminal cost \eqref{eqn_SC_NMPC_problem_formulation_terminal_cost}.
\begin{prop}
	\label{prop_terminal_cost}
The explicit value of the terminal cost is: 
	\begin{align}
			\label{eqn_SC_terminal_cost_t}
		F(\mathbf{x})  
		&= \theta^2 \left( \dfrac{\sqrt{2} \sqrt{V^2 + 2a_m r} - V}{a_m }\right)
		+	\frac{1}{3} \frac{|\theta^3|}{\omega_m} \nonumber  \\
		&+ \dfrac{\sqrt{2}\left(V^2 + 2a_m r\right)^{3/2}\left(23V^2 + 40a_m^2 + 46a_m r\right)}{240 a_m^3} \nonumber \\
		& - \left(
		\dfrac{1}{3} \dfrac{V^3}{a_m} 
		+ \dfrac{V r^2}{a_m} 
		+ \dfrac{2}{3} \dfrac{V^3 r }{a_m^2} 
		+ \dfrac{2}{15} \dfrac{V^5 }{a_m^3} 
		\right),
	\end{align} 
with $V=\sqrt{V_x^2+V_z^2} $ and $r=\sqrt{x^2+z^2}$.
	\vspace{0.5pt}
\end{prop}
\proof The proof is in Appendix \ref{appendix_terminal_cost_steps}. \hfill\ensuremath{\square}

\subsection{Stability analysis}
\label{subsec_Stability_Analysis}
For the NMPC algorithm \ref{algorithm_NMPC_SC_design}, with prediction horizon $T_p$, stage cost $L$ as in \eqref{eqn_SC_NMPC_problem_formulation_stage_cost}, terminal cost $F$ as in \eqref{eqn_SC_NMPC_problem_formulation_terminal_cost}, and the terminal set $\mathcal{X}_f$, stability can be guaranteed if these design parameters are chosen satisfying the following five stability conditions. \cite{fontes_GeneralFramework_2001}
\begin{enumerate}[label={SC\arabic*},leftmargin=22pt]
	\item \label{SC_1} The terminal set $\mathcal{X}_f$ is closed and contains the origin. 
	\item \label{SC_2} The stage cost $L$ is continuous and positive definite. 
	\item \label{SC_3} The terminal cost $F$ is positive semi-definite and continuously differentiable. 
	\item \label{SC_4} $\mathcal{X}_f$ is reachable in the prediction horizon $T_p$ from any initial states.
	\item \label{SC_Lyapunov} For all
 $t^\star \in [t_f,\infty)$ and for all $\mathbf{x}_{t} \in \mathcal{X}_f$, there exists a $ \delta >0$ such that the control function $\mathbf{u}_\text{aux}:[t^\star,t^\star+\delta] \rightarrow R^2$,  verifies $\mathbf{u}_\text{aux} \in \mathcal{U}$, makes the terminal set $\mathcal{X}_f$ invariant, and satisfies:
	\begin{equation}
		\label{eqn_SC_Lyapunov}
		\dfrac{dF(\mathbf{x}(t; t^\star,\mathbf{x}_{t},\mathbf{u}_\text{aux} ))}{dt} 
		\leq -L(\mathbf{x}(t; t^\star,\mathbf{x}_{t},\mathbf{u}_\text{aux} )).
	\end{equation}
\end{enumerate}
\begin{thm}
	\label{theorem_NMPC_scheme}
Since the model \eqref{eqn_SC_state_space_dynamics}
satisfies assumptions H$1$--H$4$ in \cite{fontes_GeneralFramework_2001}, the chosen design parameters $T_p$, $L$, $F$, and $\mathcal{X}_f$ satisfy the stability conditions \ref{SC_1}--\ref{SC_Lyapunov},
the NMPC scheme of Algorithm \ref{algorithm_NMPC_SC_design} is  stabilizing.
\end{thm}
\proof
The model clearly satisfies H$1$--H$4$ in \cite{fontes_GeneralFramework_2001} since $\mathcal{U}$ is bounded and $f(0,0)=0$.
The conditions \ref{SC_1}--\ref{SC_2} can be easily verified. 
\textcolor{black}{From the definition of $F$ in \eqref{eqn_SC_NMPC_problem_formulation_terminal_cost}, and \ref{SC_2}, $F$ is clearly positive semi-definite, and its time derivative is continuous, so \ref{SC_3} is satisfied.}
\textcolor{black}{For \ref{SC_4}, the spacecraft can rotate to the terminal set $\mathcal{X}_f$ with the maximum angle of $\pi$.}
\textcolor{black}{We provide here a sketch of the proof for the last condition \ref{SC_Lyapunov}, showing that $\dot{F}(\mathbf{x}(t; t_f,\mathbf{x}_{t_f},\mathbf{u}_\text{aux} )) = -L(\mathbf{x}(t; t_f,\mathbf{x}_{t_f},\mathbf{u}_\text{aux} ))$}.
By employing the dynamics inside the terminal set as in \eqref{eqn_SC_terminal_set_dynamics_shortened}, the time derivative of the terminal cost \eqref{eqn_SC_terminal_cost_t} is:
\begin{align}
&	\dfrac{dF(\mathbf{x}(t; t_f,\mathbf{x}_{t_f},\mathbf{u}_\text{aux} ))}{dt} 
	=
	\dfrac{\partial F(\mathbf{x}(t))}{\partial \mathbf{x}_{ts}} \cdot
	\dfrac{d\mathbf{x}_{ts}(t)}{dt} \nonumber \\
	& =
	\nabla_{\mathbf{x}_{ts}} F(\mathbf{x}(t)) \cdot
	 f_{ts}(\mathbf{x}_{ts}(t),\mathbf{u}_\text{aux}(t))
\end{align}
with the auxiliary controller as in \ref{AC_2}--\ref{AC_3}, but instead of using the time dependence, we utilize the state dependence:
\begin{align*}	\label{eqn_SC_temrinal_cost_auxiliary_controller_sign}
	&	\mathbf{u}_\text{aux}    =  [a, \omega]^\top =  \nonumber \\
	& \left\{
	\begin{array}{lll}
		[a_m, 0]^\top 										& \hspace{-6pt} \triangleq {_1}\mathbf{u}_\text{aux}& \hspace{-8pt} \text{ if }  (r,V)  \neq \bm{0}, V^2 < 2a_m r  \\
		\left[ -a_m,0 \right]^\top 							&\hspace{-6pt} \triangleq {_2}\mathbf{u}_\text{aux} & \hspace{-8pt}\text{ if }  (r,V)  \neq \bm{0}, V^2 = 2a_m r \\
		\left[ 0, -\text{sign}(\theta)w_m \right]^\top  & \hspace{-6pt} \triangleq {_3}\mathbf{u}_\text{aux} & \hspace{-8pt} \text{ if }  (r,V)  = \bm{0},  \\
	\end{array}\right.
\end{align*}
which implies \emph{Accelerate}, \emph{Decelerate}, and \emph{Rotate at origin} inside the terminal set $\mathcal{X}_f$.

The explicit equations for the partial derivative $\nabla_{\mathbf{x}_{ts}} F(\mathbf{x}(t))$ and the time derivative $\dot{F}(\mathbf{x}(t; t_f,\mathbf{x}_{t_f},\mathbf{u}_\text{aux} ))$
  when $(r,V)\neq \bm{0}$ are in Appendix \ref{appendix_terminal_cost_Lyapunov}.
%
If $(r,V)=\bm{0}$, the terminal cost \eqref{eqn_SC_terminal_cost_t} and its partial derivative are: 
\begin{equation}
F(\mathbf{x}) = \frac{1}{3} \frac{|\theta^3|}{\omega_m} , \quad 
\nabla_{\mathbf{x}_{ts}} F(\mathbf{x}(t)) = \dfrac{\text{sign}(\theta)  \theta^2}{w_m}.
\end{equation}

We now check three cases, which are equivalent to the three cases of the auxiliary controller.
\begin{enumerate}
    \item  $\mathbf{u}_\text{aux} =  {_1}\mathbf{u}_\text{aux} = 	[a_m, 0]^\top 	$,  $(r,V)  \neq \bm{0}, V^2 < 2a_m r$:
\begin{equation}
	\label{eqn_Lyapunov_proof_case_1}
	\dfrac{dF(\mathbf{x}(t; t_f,\mathbf{x}_{t_f}, {_1}\mathbf{u}_\text{aux} ))}{dt} 
	= -(r^2+\theta^2+V^2).
\end{equation}
\item  $\mathbf{u}_\text{aux} =  {_2}\mathbf{u}_\text{aux} = [-a_m, 0]^\top$, $(r,V)  \neq \bm{0}, V^2 = 2a_m r$:
\begin{equation}
		\label{eqn_Lyapunov_proof_case_2}
	\dfrac{dF(\mathbf{x}(t; t_f,\mathbf{x}_{t_f}, {_2}\mathbf{u}_\text{aux} ))}{dt} 
	= -(r^2+\theta^2+V^2).
\end{equation}
\item $\mathbf{u}_\text{aux} =  {_3}\mathbf{u}_\text{aux} = [0, -\text{sign}(\theta)\omega_m]^\top$, $(r,V)  = \bm{0}$:
\begin{equation}
		\label{eqn_Lyapunov_proof_case_3}
	\dfrac{dF(\mathbf{x}(t; t_f,\mathbf{x}_{t_f},{_3}\mathbf{u}_\text{aux} ))}{dt} 
 	=-\text{sign}^2(\theta)  \theta^2 = -\theta^2.
\end{equation}
\end{enumerate}
Eqns. \eqref{eqn_Lyapunov_proof_case_1}--\eqref{eqn_Lyapunov_proof_case_3} satisfy SC5, 
thus, completing the proof. \hfill\ensuremath{\square}

\section{SIMULATIONS}
\label{sec_simulation}
To validate the NMPC controller in Algorithm \ref{algorithm_NMPC_SC_design}, we test it on a spacecraft with the following input constraints:
\begin{equation}
	\label{eqn_SC_dynamics_input_constaint_set_sim}
	\mathcal{U} \triangleq \{\mathbf{u}\in \mathbb{R}^2 \mid -[\sqrt{2}, \pi/8]^\top \leq \mathbf{u} \leq [\sqrt{2}, \pi/8]^\top\}.
\end{equation}
The function $\theta_r(x,z) $ in \eqref{eqn_theta_ref} is \textcolor{black}{conveniently} implemented in \verb|MATLAB| or \verb|Python| with the command $\theta_r(x,z) = atan2(x,-z)$. 
\textcolor{black}{We first run a forward simulation with the controller \eqref{eqn_SC_auxiliary_controller}, then simulate NMPC controllers.}
\subsection{Initial forward simulation}
\begin{figure}[ht!] 
	\centering
	\begin{subfigure}[b]{0.49\columnwidth}
		\centering
		\includegraphics[width=0.8\textwidth]{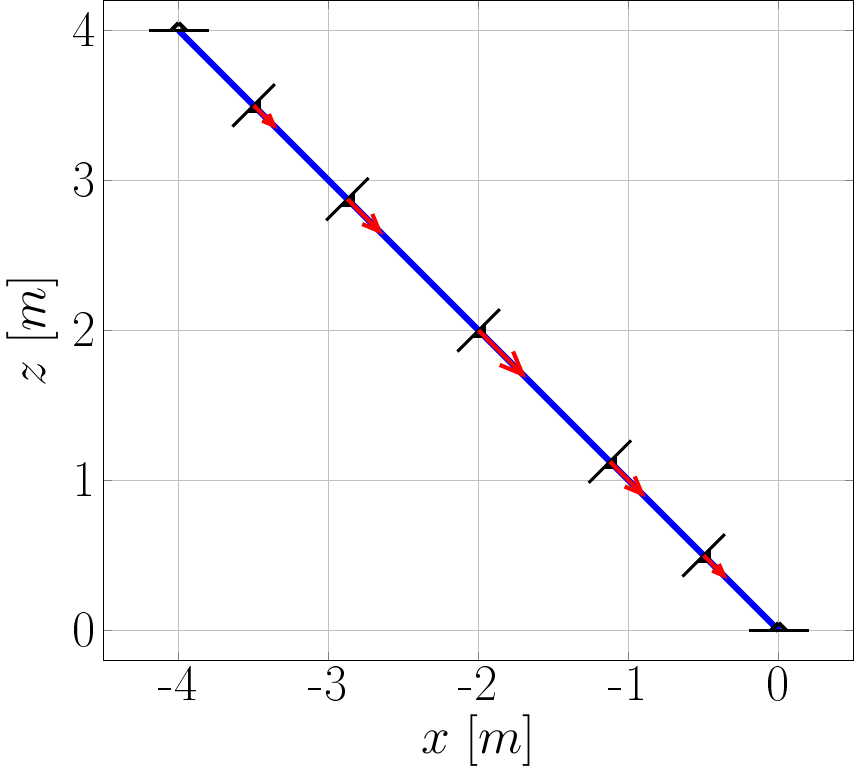}
		\caption{$x-z$}
		\label{}
	\end{subfigure}
	\hfill
	\begin{subfigure}[b]{0.49\columnwidth}
		\centering
		\includegraphics[width=\textwidth]{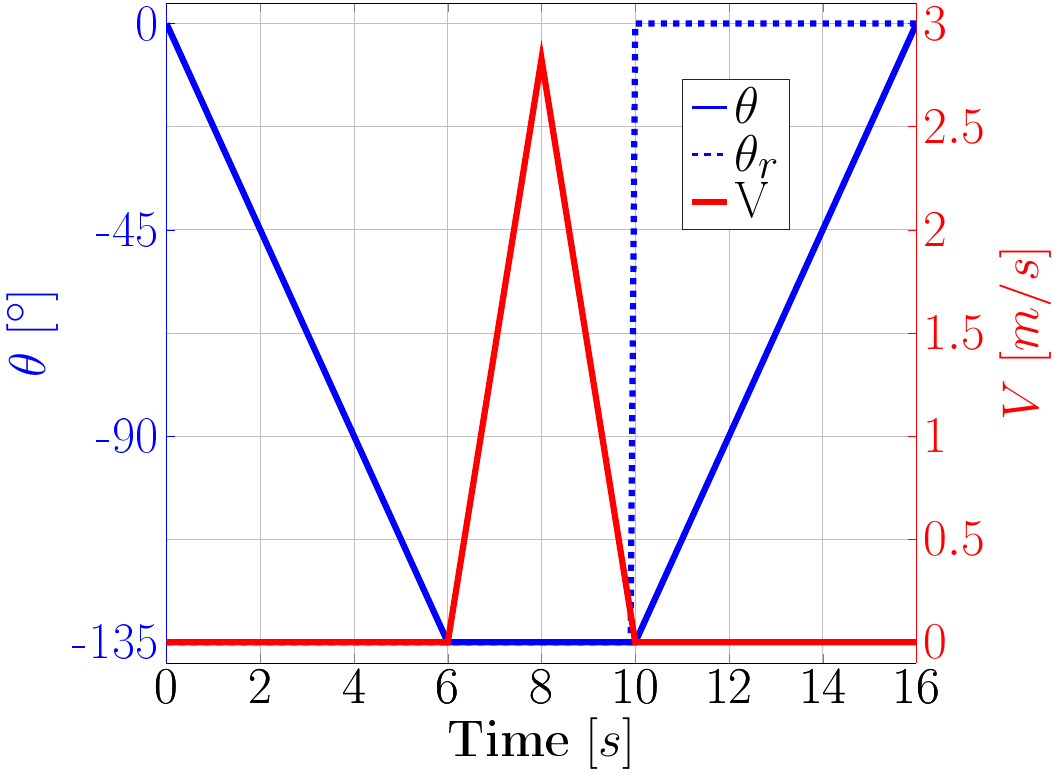}
		\caption{$\theta$ and $V$}
		\label{}
	\end{subfigure}
	\hfill
	\begin{subfigure}[b]{0.49\columnwidth}
		\centering
		\includegraphics[width=\textwidth]{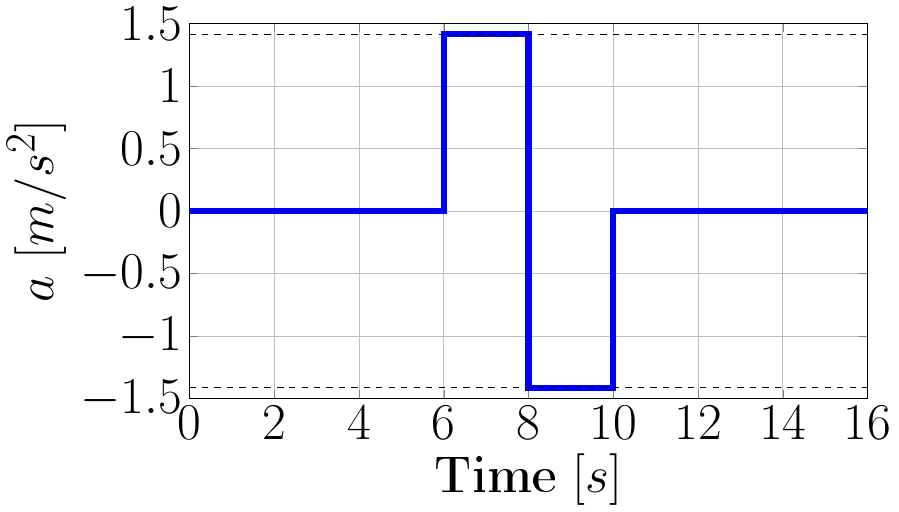}
		\caption{$a$}
		\label{}
	\end{subfigure}
	\hfill
	\begin{subfigure}[b]{0.49\columnwidth}
		\centering
		\includegraphics[width=\textwidth]{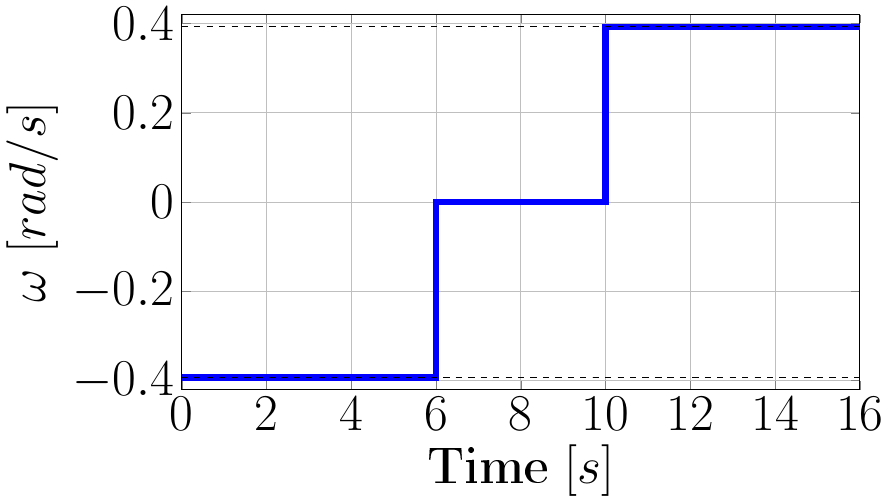}
		\caption{$\omega$}
		\label{}
	\end{subfigure}
	\caption{Simulation results for the spacecraft with the auxiliary controller, starting at $(x_0,z_0)=(-4,4)$.}
	\label{fig_SC_sim_fwd}
\end{figure}

We first simulate the motion of the spacecraft from the initial point $\mathbf{x}_A =[-4,4,0,0,0,0]^\top$ to the point $\mathbf{x} = \bm{0}$ using the auxiliary control law \eqref{eqn_SC_auxiliary_controller}. The dynamics \eqref{eqn_SC_state_space_dynamics} is discretized using the Runge-Kutta $4th$ order method with the step size $h = 0.1$, \textcolor{black}{we choose the MPC sampling time to be $\delta = 0.1s$}. The simulation results are in Fig. \ref{fig_SC_sim_fwd}.

\textcolor{black}{The reference pitch angle $\theta_r$ following \eqref{eqn_theta_ref} is $-3\pi/4 \ (rad)$, so the spacecraft takes $(-3\pi/4)/(-\pi/8)=6s$ to align $\bm{z^\mathcal{B}}$ pointing to the origin. 
Next, using \eqref{eqn_SC_T1_T2_explicit_main}, $T_1=T_2=2s$.
Finally, employing \eqref{eqn_SC_auxiliary_controller}, the rotation to $\theta=0$ at $(x,z)=(0,0)$ takes $\tfrac{\text{sign}(-(-3\pi/4)) }{\pi/8}=6s$ to complete.
}
\subsection{NMPC simulation}

For the NMPC simulation, we want to stabilize the spacecraft also from the initial state $\mathbf{x}_A =[-4,4,0,0,0,0]^\top$.
We choose prediction horizon to be $T_p = \left\lfloor \tfrac{\pi}{\omega_m \delta} \right\rfloor + 1 = 61 \text{ (steps)}$.
The OCP \eqref{eqn_SC_NMPC_problem_formulation}--\eqref{eqn_SC_NMPC_problem_formulation_constraints} is solved using the Direct Multiple Shooting technique \cite{bock_MultipleShooting_1984}, the dynamics \eqref{eqn_SC_state_space_dynamics} are discretized using the Runge-Kutta $4th$ order method with the step size $h = 0.1$.
We choose \textcolor{black}{the error tolerance} $\epsilon_r$ in Algorithm \ref{algorithm_NMPC_SC_design}, line \ref{eqn_epsilon_r} to be  $\epsilon_r=10^{-8}$.

\begin{figure}[ht!] 
	\centering
	\begin{subfigure}[b]{0.49\columnwidth}
	\centering
\includegraphics[width=\linewidth]{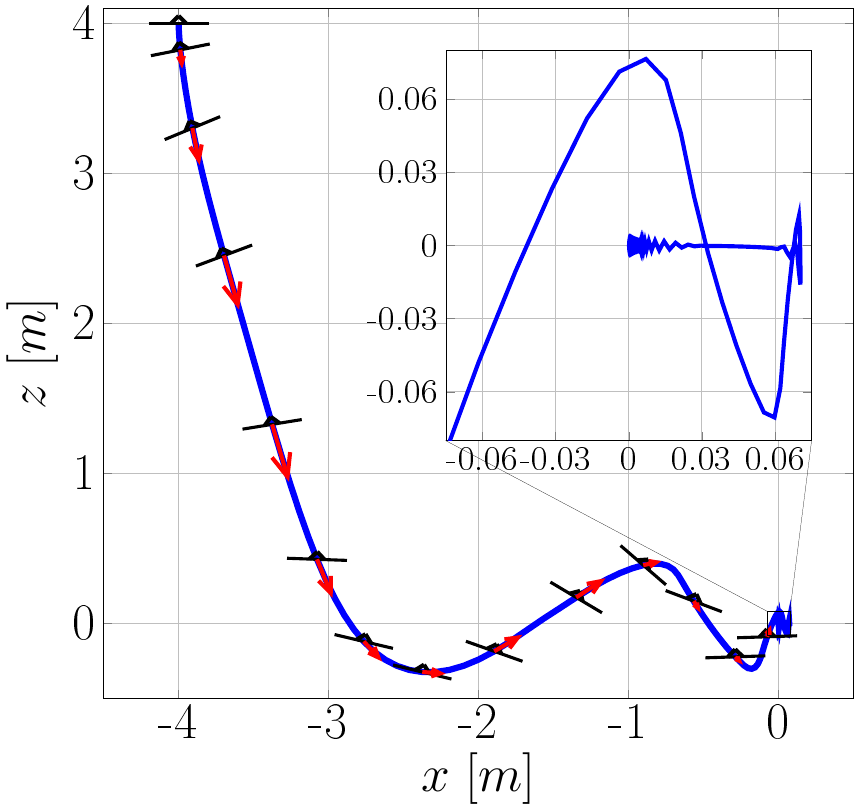}
\caption{$(x_0,z_0)=(-4,4)$}
\label{fig_SC_sim_NMPC_x_z}
	\end{subfigure}
	\hfill
	\begin{subfigure}[b]{0.49\columnwidth}
		\centering
		\includegraphics[width=0.94\linewidth]{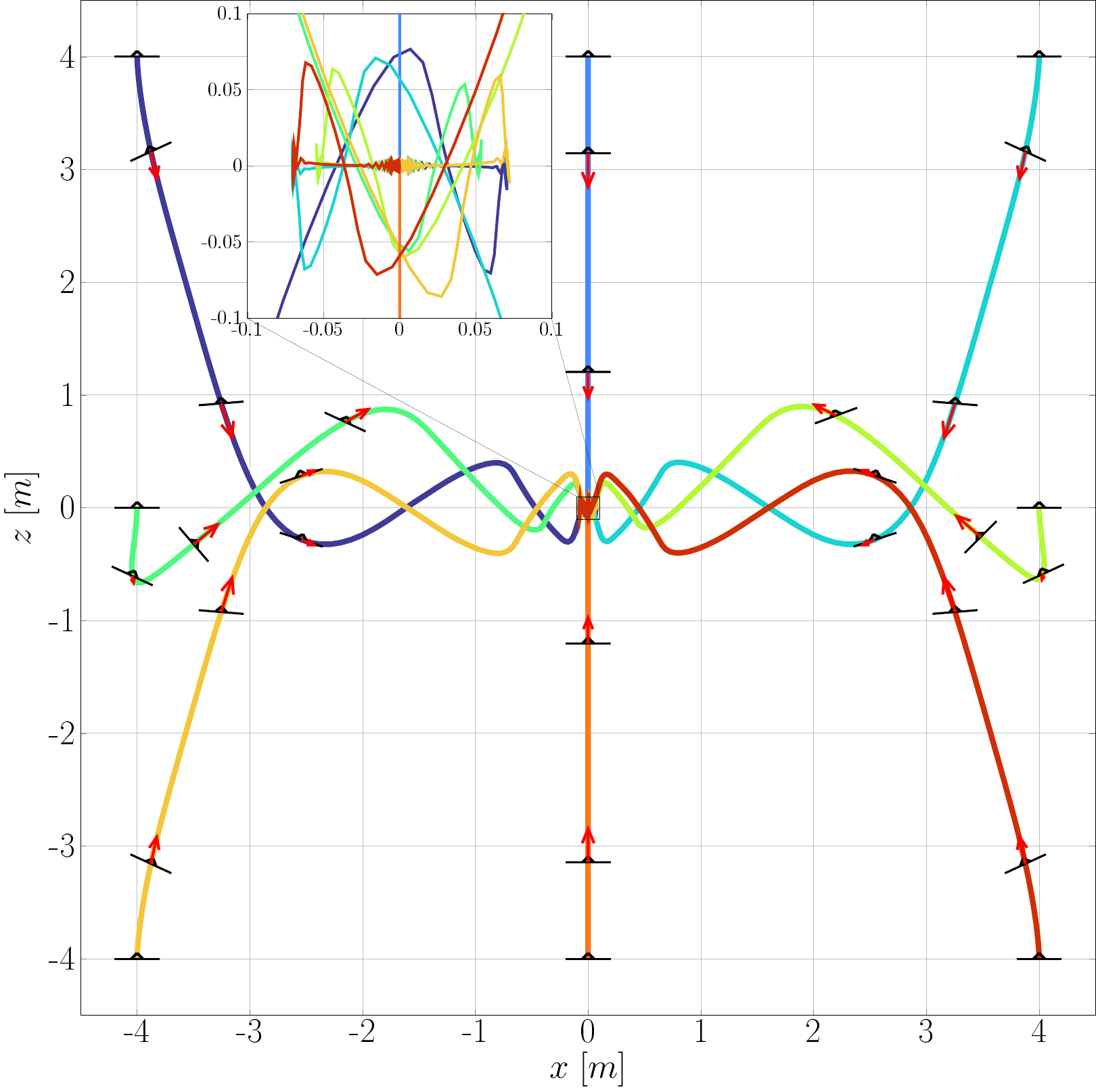}
		\caption{Various $(x_0,z_0)$}
		\label{fig_SC_sim_NMPC_x_z_bundle}
	\end{subfigure}
	\caption{Spacecraft $x-z$ trajectory.}
	\label{fig_SC_sim_NMPC_x_z_two}
\end{figure}

\begin{figure}[ht!] 
	\centering
	\begin{subfigure}[b]{0.99\columnwidth}
		\centering
		\includegraphics[width=\textwidth]{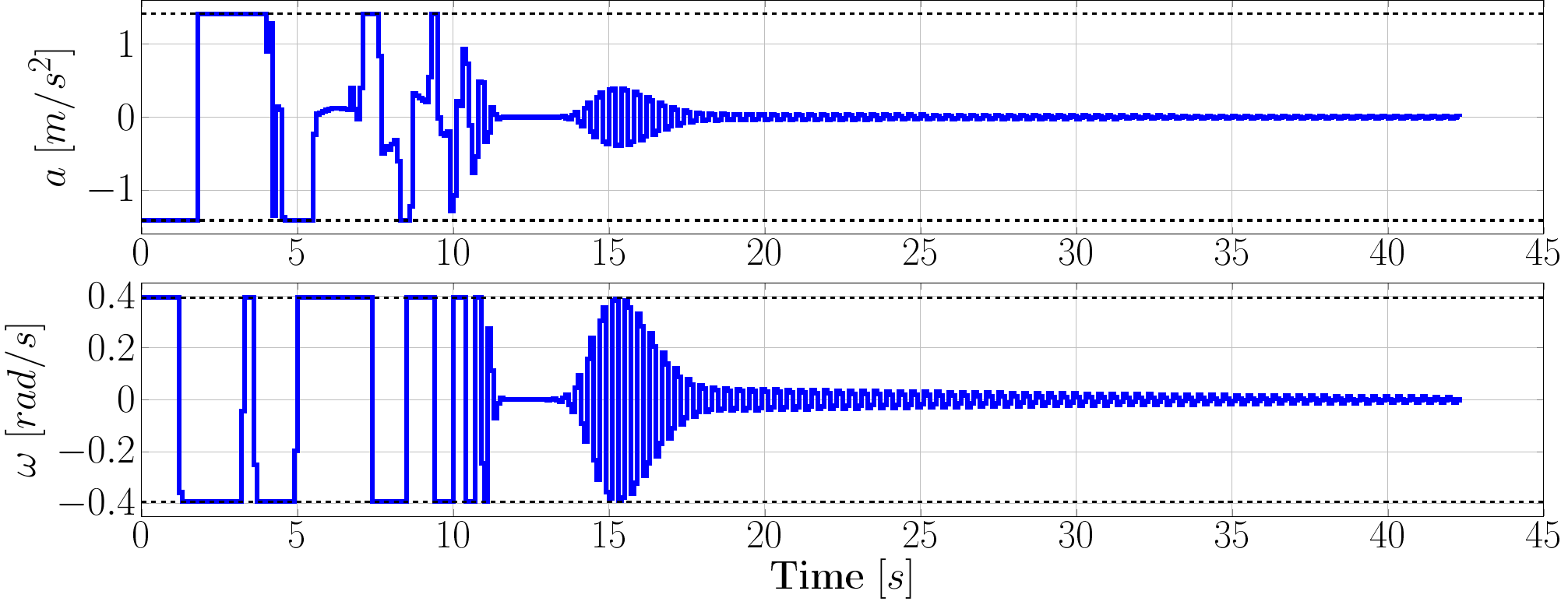}
		\caption{Control inputs $a$ and $\omega$}
		\label{fig_SC_sim_NMPC_a_omega}
	\end{subfigure}
	\hfill
	\begin{subfigure}[b]{0.49\columnwidth}
		\centering
		\includegraphics[width=\textwidth]{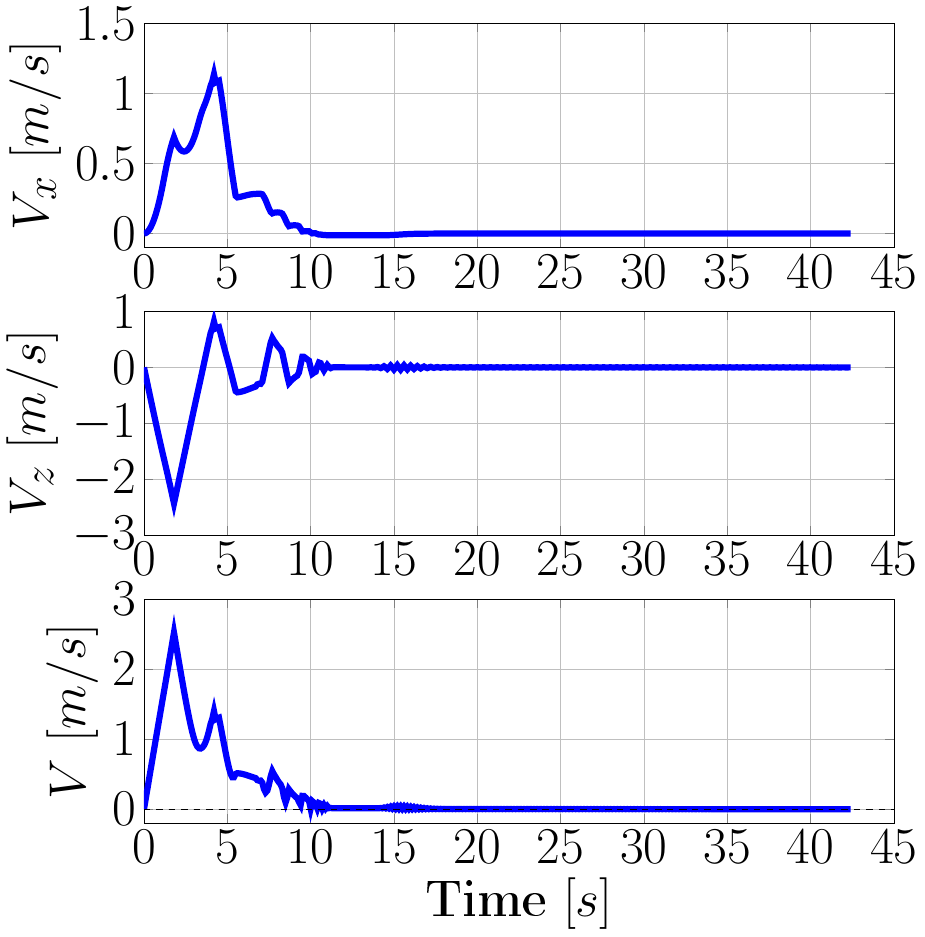}
		\caption{$V_x$, $V_z$, and $V$}
		\label{fig_SC_sim_NMPC_Vx_Vz_V}
	\end{subfigure}
	\hfill
	\begin{subfigure}[b]{0.49\columnwidth}
		\centering
		\includegraphics[width=\textwidth]{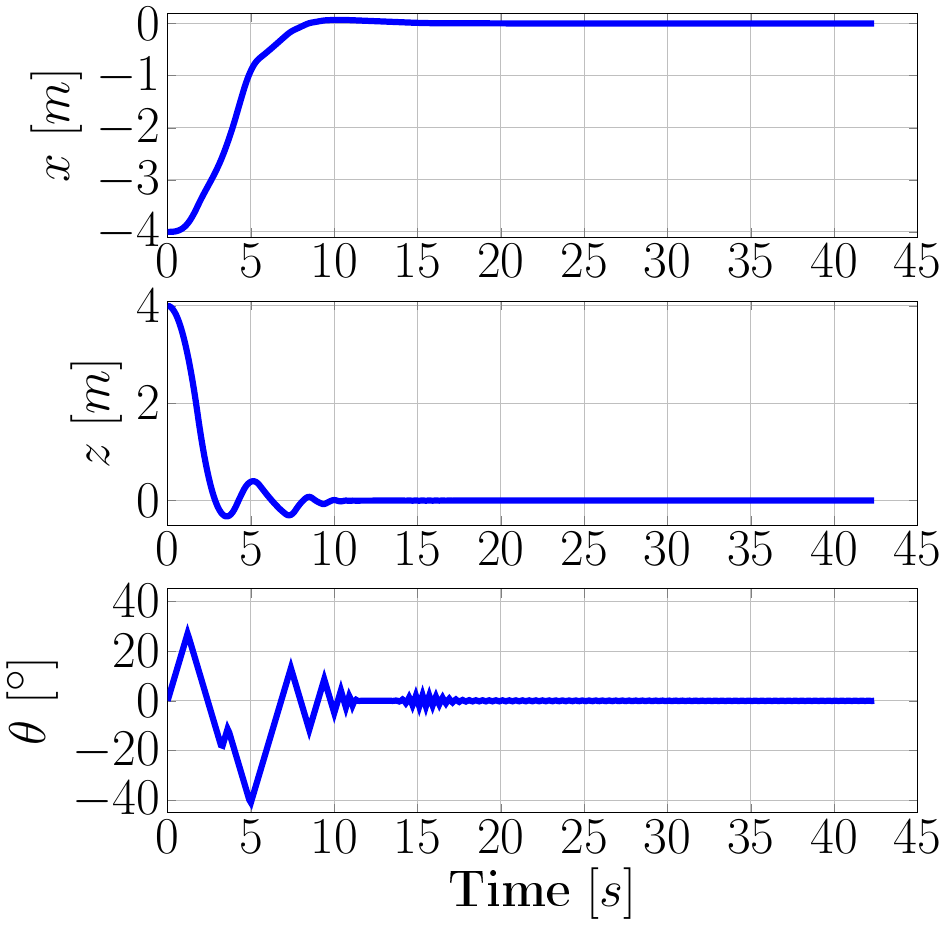}
		\caption{$x$, $z$, and $\theta$}
		\label{fig_SC_sim_NMPC_xt_zt_theta}
	\end{subfigure}
	\caption{Simulation results for the spacecraft with the NMPC controller, starting at $(x_0,z_0)=(-4,4)$.}
	\label{fig_SC_sim_NMPC}
\end{figure}

The simulation is done in MATLAB R2023a, 
using the CasADi toolbox \cite{andersson_CasADiSoftware_2019} with IPOPT solver \cite{wachter_ImplementationInteriorpoint_2006}.
\textcolor{black}{We formulate the} terminal set constraint \eqref{eqn_SC_NMPC_terminal_set_def} \textcolor{black}{on vector directions} \textcolor{black}{using the dot product and the cross product:}
\begin{subequations}
	\label{eqn_epsilon_f}
	\begin{alignat}{2}
  \bm{V_{t_f}} \times \bm{e}_z^\mathcal{B}  & = 0, \\
  \bm{r_{t_f}} \times \bm{e}_z^\mathcal{B}  & = 0, \\
 \bm{r_{t_f}} \cdot  \bm{V_{t_f}} & = -\| \bm{r_{t_f}} \| \|  \bm{V_{t_f}} \|.
\end{alignat}
\end{subequations}
\textcolor{black}{These constraints  are implemented numerically as
	$
	(1 - \epsilon_f)\|V_{t_f}\|  \leq (-V_{x_{t_f}} \sin \theta_{t_f} + V_{z_{t_f}} \cos \theta_{t_f}) \leq  \|V_{t_f}\|$,
	$- \|r_{t_f}\|               \leq (-{x_{t_f}} \sin \theta_{t_f} + {z_{t_f}} \cos \theta_{t_f})^2   \leq (-1+\epsilon_f) \|r_{t_f}\|$, 
	$- \|V_{t_f}\|\|r_{t_f}\|  	 \leq (V_{x_{t_f}} {z_{t_f}} + V_{z_{t_f}} {x_{t_f}})^2   \leq (-1+\epsilon_f) \|V_{t_f}\|\|r_{t_f}\|$, with $\epsilon_f=0$.}
We use a solution of the forward simulation, where the spacecraft is first inside the terminal set (i.e., when the pitch angle equals the reference pitch angle) as the initial guess for the first NMPC iteration. Other NMPC iterations use the warm-start method from the solution of the previous iteration.  
\textcolor{black}{We choose the sampling time $\delta = 0.1s$.}
The simulation results are in Figs. \ref{fig_SC_sim_NMPC_x_z}--\ref{fig_SC_sim_NMPC}. 

\vspace{-10pt}
The control inputs $a$ and $\omega$ are in Fig. \ref{fig_SC_sim_NMPC_a_omega}, satisfying the input constraints \eqref{eqn_SC_dynamics_input_constaint_set_sim}.
In Figure \ref{fig_SC_sim_NMPC_x_z}, we see that the spacecraft approaches the origin from the right, with the pitch angle $\theta$ chattering around $0^\circ$ (Fig. \ref{fig_SC_sim_NMPC_xt_zt_theta}). Since the spacecraft has no lateral acceleration, it executes damped oscillations around its center of mass by attenuating the angular rate $\omega$ (Fig. \ref{fig_SC_sim_NMPC_a_omega}). 
All the five states of the spacecraft are converging to $0$, which illustrates the stability of the proposed NMPC scheme.

Figure \ref{fig_SC_sim_NMPC_x_z_bundle} shows the $x-z$ trajectory of the spacecraft for $(x_0,z_0) \in \{-4,0,4\} \times \{-4,0,4\} \setminus (0,0)$. \textcolor{black}{An interesting observation is the spacecraft executes the final damping maneuvers on the half plane $x \geq 0$ or  $x \leq 0$ that is the opposite of the other half-plane in which it started at $x_0$.}
\section{CONCLUSIONS AND FUTURE WORKS}
\label{sec_conclusion}
This paper presents an NMPC scheme for a class of drift nonholonomic systems with lateral acceleration restrained, where the terminal set is a region of a subspace, which contains the modes where the nonholonomic system is easily controllable.
For future work, we will investigate how  to implement a similar technique to amulticopter, which is another drifting nonholonomic system.
\bibliographystyle{IEEEtran}
\bibliography{CDC2024}   

\clearpage
  \appendices
\section{Steps to calculate the terminal cost}
\label{appendix_terminal_cost_steps}

We decompose the terminal cost \eqref{eqn_SC_NMPC_problem_formulation_terminal_cost} into four terms, which are equivalent to the time intervals $\tau \in [t_f,t_f+t_1]$, $\tau \in [t_f+t_1,t_f+t_2]$, $\tau \in [t_f+t_2,t_f+t_3]$, and $\tau \in [t_f+t_3,\infty]$, and we respectively denote each term as $F_1$, $F_2$, $F_3$, and $F_4$:
\begin{align*}
	\label{eqn_SC_NMPC_problem_formulation_terminal_cost_split}
	F((\bar{\mathbf{x}}(t_f))   
	&= \underbrace{\int_{t_f}^{t_f+t_1} L((\bar{\mathbf{x}}(\tau))   d\tau }_{F_1}
	+ \underbrace{\int_{t_f+t_1}^{t_f+t_2} L((\bar{\mathbf{x}}(\tau))   d\tau }_{F_2} \nonumber \\
	&+ \underbrace{\int_{t_f+t_2}^{t_f+t_3} L((\bar{\mathbf{x}}(\tau))   d\tau }_{F_3}
	+ \underbrace{\int_{t_f+t_3}^{\infty} L((\bar{\mathbf{x}}(\tau))   d\tau }_{F_4}.
\end{align*}

We now proceed to calculate the explicit values of $F_1$, $F_2$, $F_3$, $F_4$, and the terminal cost $F$.

The velocity at $t_f+t_1$, which we call the maximum velocity $V_m$, is: 
\begin{equation*}
	\label{eqn_SC_terminal_cost_Vm}
	V_m = V_{t_f} + \int_{t_f}^{t_f+t_1} a(\tau) d\tau = V_{t_f} + a_m t_1. 
\end{equation*}

From $t_f+t_1$ to $t_f+t_2$, by looking at the Fig. \ref{fig_SC_terminal_set_acc_vel}, we can obtain the relation between $t_2$ and $t_1$ as follows
\begin{equation}
	\label{eqn_SC_terminal_cost_t2_t1_relation}
	t_2 - t_1 = \dfrac{V_m}{a_m} =  \dfrac{V_{t_f} + a_m t_1 }{a_m}
	\Rightarrow 
	t_2 = 2 t_1 + \dfrac{V_{t_f}}{a_m}.
\end{equation}

We want at $t_f+t_2$, the distance $r(t_f+t_2)$ is $0$, which means the area of below the function $V(t)$ over $t$ as in Fig. \ref{fig_SC_terminal_set_acc_vel}, equals to $r_{t_f}$:
\begin{equation}
	\label{eqn_SC_terminal_cost_r_area}
	\begin{aligned}[b]
		r_{t_f} =& \ \dfrac{1}{2} (V_{t_f} + V_m) t_1 + \dfrac{1}{2} V_m (t_2 - t_1) \\
		=& \ \dfrac{1}{2} \left(4V_{t_f}t_1 + 2a_m t_1^2 +\dfrac{V_{t_f}^2}{a_m} \right) \\
		=&\ 2 V_{t_f}t_1 + a_m t_1^2 + \dfrac{V_{t_f}^2}{2a_m}.
	\end{aligned}
\end{equation}

Rewriting the terms of \eqref{eqn_SC_terminal_cost_r_area} and solving for $t_1$, we obtain:
\begin{equation}
	\label{eqn_SC_terminal_cost_t1_eqn}
	a_m t_1^2 + 2 V_{t_f}t_1 + \left(  \dfrac{V_{t_f}^2}{2a_m} - r_{t_f} \right) = 0,
\end{equation}
with $\Delta = (2 V_{t_f})^2 - 4 (a_m) ( \tfrac{V_{t_f}^2}{2a_m} - r_{t_f})
	= 2 (V_{t_f}^2 + 2a_m r_{t_f}).$
Since $V_{t_f} \geq 0$, $a_m \geq 0$, and $r_{t_f} \geq 0$, which lead to $\Delta \geq 0$, hence \eqref{eqn_SC_terminal_cost_t1_eqn} always has at least one solution.

By solving \eqref{eqn_SC_terminal_cost_t1_eqn}, combing with \eqref{eqn_SC_terminal_cost_t2_t1_relation}, the explicit formulas for $t_1$ and $t_2$ are: 
\begin{subequations}
	\label{eqn_SC_t1_t2_explicit}
	\begin{align}
		\label{eqn_SC_t1_explicit}
		t_1 &= -\dfrac{V_{t_f}}{a_m} + \dfrac{\sqrt{V_{t_f}^2 + 2a_m r_{t_f}}}{a_m \sqrt{2}}, \\	
		\label{eqn_SC_t2_explicit}
		t_2 &= -\dfrac{V_{t_f}}{a_m} + \dfrac{\sqrt{2} \sqrt{V_{t_f}^2 + 2a_m r_{t_f}}}{a_m }.
	\end{align}
\end{subequations}

\begin{rem}
Note that, to make the auxiliary controller valid, i.e. we want the spacecraft to accelerate then decelerate \emph{inside the terminal set $\mathcal{X}_f$}, we must have the condition $t_1 \geq 0$, which is equivalent to the condition: 
\begin{equation}
	V_{x_{t_f}}^2 +V_{z_{t_f}}^2 \leq 2 a_m \sqrt{x_{t_f}^2+z_{t_f}^2}.
\end{equation}
\end{rem}

For $t \in [t_f, t_f + t_1]$, with the acceleration $a(t) = a_m$ and the pitch angle $\theta(t) = \theta_{t_f} = \theta_r$, the velocity and the distance of the spacecraft to the origin are:
	\begin{align*}
		V(t) 
		&= V_{t_f} +  \int_{t_f}^{ t}   a(\tau)   d\tau
		= V_{t_f} +  \int_{t_f}^{t}  ( a_m )  d\tau \nonumber  \\
		&= V_{t_f} +  a_m (t-t_f),
	\end{align*}
	\begin{align*}
		r(t) =&\ r_{t_f} - \int_{t_f}^{ t}   V(\tau)   d\tau
		=  r_{t_f} - \int_{t_f}^{ t}  [V_{t_f} + a_m (\tau -t_f)]  d\tau \nonumber \\
		=&\ r_{t_f} - V_{t_f}(t-t_f) -  \dfrac{a_m}{2} (t -t_f)^2. 
	\end{align*}
The first terminal cost term $F_1$ is now calculated as
\begin{equation*}
	\label{eqn_SC_terminal_cost_F1_tf_tf_t1}
	\begin{aligned}[b]
		F_1 =&\   \int_{t_f}^{t_f+t_1}  (r(\tau)^2 +  V(\tau)^2+ \theta(\tau)^2 )   d\tau \\
		=&\ \theta_{t_f}^2 t_1 
		+ \hspace{-2pt} \int_{t_f}^{t_f+t_1}  \hspace{-3pt}
		\left\{
		\left[ r_{t_f}  - V_{t_f}(\tau-t_f) -  \dfrac{a_m}{2} (\tau -t_f)^2 \right]^2 \right.\\
		&+ \left.  \left[ V_{t_f} + a_m (\tau-t_f) \right]^2 
		\right\}   d\tau \\
		=&\ (\theta_{t_f}^2  + r_{t_f}^2  +  V_{t_f}^2) t_1
		- V_{t_f}(r_{t_f}-a_m)t_1^2 \\
		&+ \frac{1}{3} (V_{t_f}^2 + a_m^2 -  a_m r_{t_f})  t_1^3
		+ \frac{1}{4} a_m V_{t_f} t_1^4
		+  \frac{1}{20} a_m^2 t_1^5.
	\end{aligned}
\end{equation*}

At $t = t_f + t_1$, the velocity and the distance are 
\begin{subequations}
	\label{eqn_SC_terminal_cost_V_r_tf_t1}
	\begin{align*}
		V( t_f + t_1 ) &= V_{t_f} + a_m t_1, \\
		r(t_f + t_1) &=   r_{t_f}  - V_{t_f}t_1 - \dfrac{a_m}{2} t_1^2.
	\end{align*}
\end{subequations}

For $t \in [t_f+t_1, t_f + t_2]$, with the acceleration $a(t) = -a_m$ and the pitch angle $\theta(t) = \theta_{t_f} = \theta_r$, the velocity and the distance of the spacecraft to the origin are:
\begin{equation*}
	\label{eqn_SC_terminal_cost_V_tf_t1_tf_t2}
	\begin{aligned}[b]
		V(t) 
		&= V(t_f+t_1) +  \int_{t_f+t_1}^{ t}   a(\tau)   d\tau \\
		&=  V_{t_f} + a_m t_1  -a_m (t-t_f-t_1)  \\
		&= V_{t_f} -a_m (t-t_f-2t_1),
	\end{aligned}
\end{equation*}
\begin{equation*}
	\label{eqn_SC_terminal_cost_r_tf_t1_tf_t2}
	\begin{aligned}[b]
		r(t) =&\ r(t_f+t_1) - \int_{t_f+t_1}^{ t}   V(\tau)   d\tau \\
		=&\  r(t_f+t_1)
		- \int_{t_f+t_1}^{ t}  [V_{t_f} -a_m (\tau -t_f-2t_1)]  d\tau \\
		=&\ \left( r_{t_f} - V_{t_f}t_1 -  \dfrac{a_m}{2} t_1^2  \right)  
		- V_{t_f}(t-t_f-t_1) \\
		& + \dfrac{a_m}{2}  (t-t_f-2t_1)^2 - \dfrac{a_m}{2} t_1^2 \\
		= &\ r_{t_f} -  V_{t_f} (t-t_f) - a_m t_1^2 
		+ \dfrac{a_m}{2}  (t-t_f-2t_1)^2.
	\end{aligned}
\end{equation*}
The second terminal cost term $F_2$ is in \eqref{eqn_SC_terminal_cost_F2_tf_t1_tf_t2}.
\begin{figure*}[t]
	\begin{align}
			\label{eqn_SC_terminal_cost_F2_tf_t1_tf_t2}
		F_2  =&\  \int_{t_f+t_1}^{t_f+t_2}  (r(\tau)^2 +  V(\tau)^2+ \theta(\tau)^2 )   d\tau  \nonumber \\
		=&\ \theta_{t_f}^2 (t_2-t_1) +  \int_{t_f+t_1}^{t_f+t_2}  
		\left\{
		\left[ r_{t_f} -  V_{t_f} (\tau-t_f) - a_m t_1^2 
		 + \dfrac{a_m}{2}  (\tau-t_f-2t_1)^2  \right]^2
		+  
		\left[V_{t_f} -a_m (\tau-t_f-2t_1) \right]^2 
		\right\}   d\tau  \nonumber  \\
		=&\ \theta_{t_f}^2 (t_2-t_1) 
		+ \frac{1}{3} \left(t_2 - t_1 \right) 
		\left\{  \left[\frac{13 t_1^4}{20}-\frac{47 t_2 t_1^3}{20}+\left(\frac{73 t_2^2}{20}+7\right) t_1^2+\left(-\frac{27}{20} t_2^3-5 t_2\right) t_1+\frac{3 t_2^4}{20}+t_2^2\right] a_m^{2} \right.  \nonumber \\
		& +\left[\frac{V_{t_f} t_1^{3}}{4}+\left[\frac{V_{t_f} t_2}{4}+r_{t_f}\right] t_1^{2}+\left[\frac{13}{4} V_{t_f} t_2^{2}-5 r_{t_f} t_2+9 V_{t_f}\right] t_1-\frac{3 V_{t_f} t_2^{3}}{4}+r_{t_f} t_2^{2}-3 V_{t_f} t_2\right] a_m   \nonumber \\
		& \Biggl. +V_{t_f}^{2} t_1^{2}+V_{t_f} \left[V_{t_f} t_2-3 r_{t_f}\right] t_1+V_{t_f}^{2} t_2^{2}-3 V_{t_f} r_{t_f} t_2+3 V_{t_f}^{2}+3 r_{t_f}^2  \Biggl\} . 
	\end{align}
\end{figure*}

For $t \in [t_f + t_2, t_f + t_3]$, the spacecraft rotates back to $\theta = 0$ from $\theta_{t_f}$, where the pitch angle and the needed for the rotation are:
\begin{subequations}
	\begin{align}
		\theta(t) & = \theta_{t_f} - \text{sign}(\theta_{t_f}) \omega_m(t - t_f -t_2), \\
		t_3 - t_2 & =  \frac{ |\theta_{t_f}|}{ \omega_m}.
	\end{align}
\end{subequations}
Therefore, the terminal cost term for the rotation back to $\theta=0$ from $\theta=\theta_{t_f}$ at the origin $(x,z)=(0,0)$ is:
\begin{equation}
	\label{eqn_SC_terminal_cost_F3_tf_t2_tf_t3}
	\begin{aligned}[b]
		F_3 = &\
		\int_{t_f + t_2}^{t_f + t_3} [ \theta_{t_f}^2 - 2 \theta_{t_f}\text{sign}(\theta_{t_f}) \omega_m(t - t_f -t_2) \\
		& + (\text{sign}(\theta_{t_f})\omega_m(t - t_f -t_2))^2 ] dt 
		= 	\frac{1}{3} \frac{|\theta_{t_f}^3|}{\omega_m}. 
	\end{aligned}
\end{equation}

When $t\in [t+t_4, \infty]$, all the states are zeros (i.e. $\mathbf{x}=\bm{0}$), therefore, 
\begin{equation}
	F_4 = 0.
\end{equation}

Finally, the terminal cost is presented in \eqref{eqn_SC_terminal_cost}, which completes the proof.
	\begin{align}
	\label{eqn_SC_terminal_cost}
		&	F(\mathbf{x}_{t_f})  
		= F_1+F_2+F_3 + F_4  \nonumber \\
		&	= \theta_{t_f}^2 \left(-\dfrac{V_{t_f}}{a_m} + \dfrac{\sqrt{2} \sqrt{V_{t_f}^2 + 2a_m r_{t_f}}}{a_m }\right) 		
		+	\frac{1}{3} \frac{|\theta_{t_f}^3|}{\omega_m} \nonumber \\
		&\	+ \frac{\sqrt{2}\left(V_{t_f}^2 + 2a_m r_{t_f}\right)^{3/2}\left(23V_{t_f}^2 + 40a_m^2 + 46a_m r_{t_f}\right)}{240 a_m^3} \nonumber \\
		&\ - \left(
		\dfrac{1}{3} \dfrac{V_{t_f}^3}{a_m} 
		+ \dfrac{V_{t_f} r_{t_f}^2}{a_m} 
		+ \dfrac{2}{3} \dfrac{V_{t_f}^3 r_{t_f} }{a_m^2} 
		+ \dfrac{2}{15} \dfrac{V_{t_f}^5 }{a_m^3} 
		\right).
	\end{align}

\section{Explicit time derivative of the terminal cost}
\label{appendix_terminal_cost_Lyapunov}
The explicit equation for the partial derivative $\nabla_{\mathbf{x}_{ts}} F(\mathbf{x}(t))$ and the time derivative $\dot{F}(\mathbf{x}(t; t_f,\mathbf{x}_{t_f},\mathbf{u}_\text{aux} ))$
when $(r,V)\neq \bm{0}$ are in \eqref{eqn_terminal_cost_partial_derivative} and \eqref{eqn_terminal_cost_time_derivative}.

\begin{figure*}[t]
	If $(r,V)\neq \bm{0}$: \\
	\begin{equation}
		\label{eqn_terminal_cost_partial_derivative}
		\nabla_{\mathbf{x}_{ts}} F(\mathbf{x}(t))
		=
		\begin{bmatrix}
			\textcolor{black}{\dfrac{23 \sqrt{2} \sigma_2^{3 / 2}}{120 a_m^2}-\dfrac{2 V^3}{3 a_m^2}+\dfrac{\sqrt{2} \theta^2}{\sqrt{\sigma_2}}-\dfrac{2 V r}{a_m}+\dfrac{\sqrt{2} \sqrt{\sigma_2} \sigma_1}{80 a_m^2}} \\
			\textcolor{black}{\dfrac{\text{sign}(\theta)  \theta^2}{w_m}-2 \theta\left(\dfrac{V}{a_m}-\dfrac{\sqrt{2} \sqrt{\sigma_2}}{a_m}\right)}\\
			\textcolor{mygreen}{\dfrac{23 \sqrt{2} V \sigma_2^{3 / 2}}{120 a_m^3}-\dfrac{V^2}{a_m}-\dfrac{2 V^4}{3 a_m^3}-\dfrac{r^2}{a_m}-\dfrac{2 V^2 r}{a_m^2}-\theta^2\left(\dfrac{\sqrt{\sigma_2} -\sqrt{2} V}{a_m \sqrt{\sigma_2}}\right)+\dfrac{\sqrt{2} V \sqrt{\sigma_2} \sigma_1}{80 a_m^3}}
		\end{bmatrix},
	\end{equation}
	\begin{equation*}
		\text{with } \sigma_1=23 V^2+40 a_m^2+46 r a_m, \text{ and } \sigma_2=V^2+2 a_m r.
	\end{equation*}
\end{figure*}
\begin{figure*}[t]
	\begin{align}
		\label{eqn_terminal_cost_time_derivative}
		\dfrac{dF(\mathbf{x}(t; t_f,\mathbf{x}_{t_f},\mathbf{u}_\text{aux} ))}{dt} 
		=& \
		\nabla_{\mathbf{x}_{ts}} F(\mathbf{x}(t)) \cdot
		f_{ts}(\mathbf{x}_{ts}(t),\mathbf{u}_\text{aux}(t))  \nonumber \\
		= &\
		-V \left(
		\textcolor{black}{\dfrac{23 \sqrt{2} \sigma_2^{3 / 2}}{120 a_m^2}-\dfrac{2 V^3}{3 a_m^2}+\dfrac{\sqrt{2} \theta^2}{\sqrt{\sigma_2}}-\dfrac{2 V r}{a_m}+\dfrac{\sqrt{2} \sqrt{\sigma_2} \sigma_1}{80 a_m^2}} 
		\right) \nonumber \\
		& +\ \omega 
		\left[
		\textcolor{black}{\dfrac{\text{sign}(\theta)  \theta^2}{w_m}-2 \theta\left(\dfrac{V}{a_m}-\dfrac{\sqrt{2} \sqrt{\sigma_2}}{a_m}\right)}
		\right] \nonumber \\
		&	+\  a \left[
		\textcolor{mygreen}{\dfrac{23 \sqrt{2} V \sigma_2^{3 / 2}}{120 a_m^3}-\dfrac{V^2}{a_m}-\dfrac{2 V^4}{3 a_m^3}-\dfrac{r^2}{a_m}-\dfrac{2 V^2 r}{a_m^2}-\theta^2\left(\dfrac{\sqrt{\sigma_2} -\sqrt{2} V}{a_m \sqrt{\sigma_2}}\right)+\dfrac{\sqrt{2} V \sqrt{\sigma_2} \sigma_1}{80 a_m^3}}
		\right].
	\end{align}
\end{figure*}

\section{Auxiliary controller analysis when the spacecraft pointing away from the origin}
\label{appendix_second_auxiliary_controller}
This appendix presents briefly the results when we use the stabilizing strategy SS' where the vertical axis $\bm{z^\mathcal{B}}$ points away from (rather than pointing to) the origin and then the spacecraft moves backwards to the origin.

The second stabilizing strategy SS' is defined as follows: 
\begin{enumerate}[label={SS'\arabic*},leftmargin=28pt]
	\item \label{eqn_SS_comma_1} Rotate the spacecraft so that its vertical axis $\bm{z^\mathcal{B}}$ points \emph{away from} the origin .	
	\item Move \emph{backwards} the target with an acceleration profile $a'_\text{ref}(t)$ until reaching the origin of the plane $(x,z) = (0,0)$.
	\item Rotate back to  $\theta = 0$ and then stop.
\end{enumerate}
With this strategy, the reference pitch is $\theta'_r(x,z) = \theta_r(x,z) - \pi$, and the second auxiliary controller is: 
\begin{align}
	\label{eqn_SC_2nd_auxiliary_controller}
	&	\mathbf{u}'_\text{aux}    = [a, \omega]^\top = \nonumber  \\
	&	\left\{
	\begin{array}{lll}
		a = 0,                 &\omega = -\omega_m &\text{ if } \theta  > \theta'_r;    \\
		a = 0,                 &\omega = \omega_m  &\text{ if } \theta  < \theta'_r;  \\
		a = a'_{r}(t), &\omega =0          &\text{ if }  \theta = \theta'_r, (r,V)\neq \bm{0};    \\
		a=0,                   &\omega =0          &\text{ if } \mathbf{x}=\bm{0}; \\
	\end{array}\right.
\end{align}
with the acceleration profile:
\begin{align*}
	a'_r(t) =  
	\left\{
	\begin{array}{lll}
		-a_m  & \hspace{-5pt} \text{if } t \in [T_0, T_0 + T_1]  &\text{(Decelerate)} \\
		a_m & \hspace{-5pt}  \text{if } t \in [T_0 + T_1, T_0 + T_2]   &\text{(Accelerate)} 
	\end{array}\right.
\end{align*}

 The auxiliary strategy SS' is equivalently defined as:
\begin{enumerate}[label={AC'\arabic*},leftmargin=29pt]
	\item Rotate to $\mathcal{X}_{f}$, i.e. rotate to align the pitch angle $\theta = \theta'_r(x,z)$, with the angular rate $\omega = \text{sign}(\theta'_r-\theta) w_m$.
	\item \label{AC_comma_2} Reverse straight to $(x,z)=(0,0)$ within $\mathcal{X}_f $ following the acceleration profile: $a=-a_m \text{ for } t \in [t_f,t_f+t_1]$ and $a= a_m \text{ for } t \in [t_f+t_1,t_f+t_2]$.
	\item \label{AC_comma_3} Rotate to $(x,z,\theta)=(0,0,0)$ with the angular rate  $\omega = -\text{sign}(\theta) w_m$ for $t \in [t_f+t_2,t_f+t_3]$ and stop.
\end{enumerate}

\begin{figure}[H]
	\centering
	\includegraphics[width=0.4\linewidth]{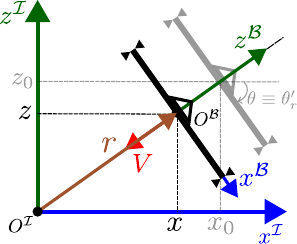}
	\caption{The spacecraft when it is inside the terminal set $\mathcal{X}_f$ and points away from origin.}
	\label{fig_Spacecraft_2D_Full_Model_r_V_reverse}
\end{figure}
Figure \ref{fig_Spacecraft_2D_Full_Model_r_V_reverse} shows a snapshot of the spacecraft when it ``reverses'' inside the terminal set \textcolor{black}{backward the origin}, with the velocity vector $\bm{V}$ in red, and the displacement vector $\bm{r}$ in brown. Within the terminal set $\mathcal{X}_f$, we obtain simpler dynamics for the spacecraft \textcolor{black}{defined in \eqref{eqn_SC_state_space_dynamics}}:
\begin{subequations}
	\label{eqn_SC_terminal_set_dynamics_reverse}
	\begin{align}
		\dot{r}(t) &= V(t),  \\
		\dot{\theta}(t) &= \omega(t), \\
		\dot{V}(t) &=  a(t), 
	\end{align}
\end{subequations}
\begin{figure}[H]
	\centering
	\begin{subfigure}[b]{0.49\linewidth}
		\centering
		\includegraphics[width=\textwidth]{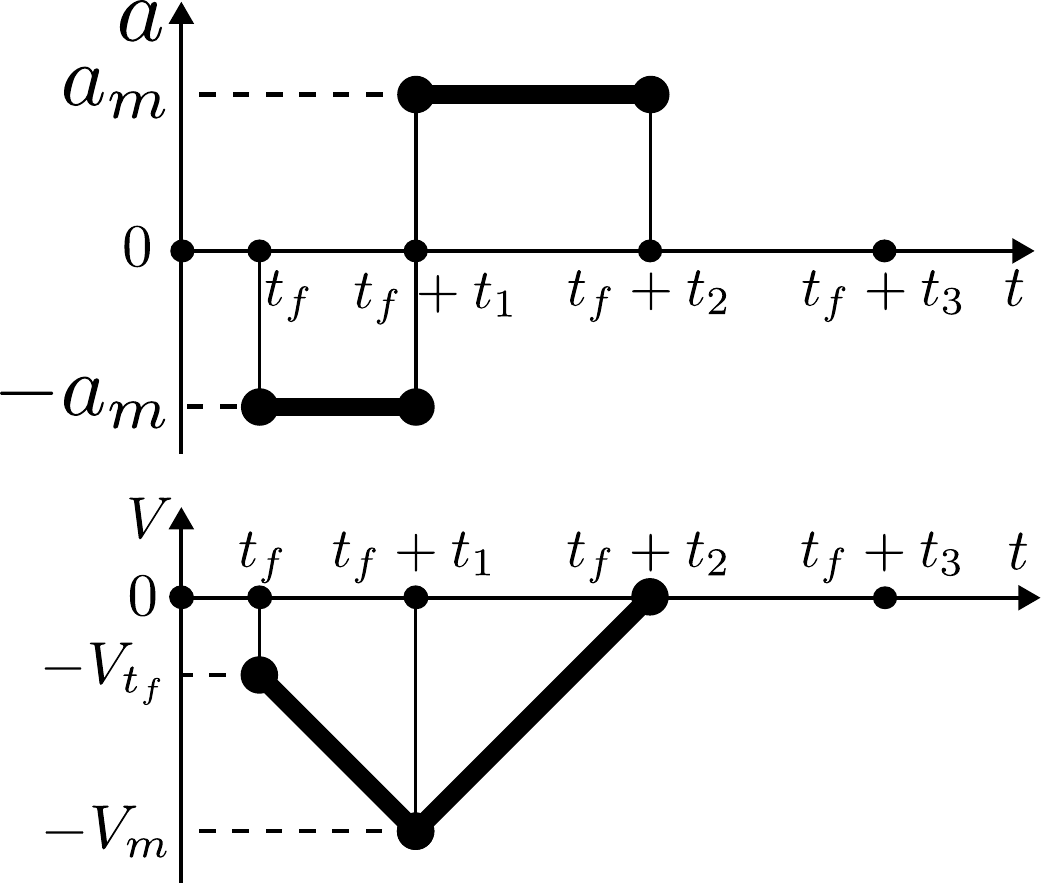}
		\caption{$a$ and $V$}
		\label{fig_SC_terminal_set_acc_vel_reverse}
	\end{subfigure}
	\hfill
	\begin{subfigure}[b]{0.49\linewidth}
		\centering
		\includegraphics[width=\textwidth]{./pics/Spacecraft/Spacecraft_terminal_set_omega_theta_profile_short_1_case.pdf}
		\caption{$\omega$ and $\theta$}
		\label{fig_SC_terminal_set_omega_theta_reverse}
	\end{subfigure}
	\caption{The auxiliary stabilizing controller AC' and states inside the terminal set.}
	\label{fig_SC_terminal_set_acc_vel_omega_theta_reverse}
\end{figure}
The difference lies in Fig. \ref{fig_SC_terminal_set_acc_vel_reverse}, where \textcolor{black}{algebraic value of the velocity vector} $V$ is negative since the spacecraft is reversing in its body frame.

The Lemma \ref{lemma_t1_t2_t3} and Proposition \ref{prop_terminal_cost} also hold with this auxiliary strategy, therefore, the Theorem \ref{theorem_NMPC_scheme} holds.

\end{document}